\documentclass[review]{elsarticle}

\usepackage{lineno}
\usepackage[utf8x]{inputenc}
\usepackage[T1]{fontenc}
\usepackage[figuresright]{rotating}
\usepackage{xcolor}
\usepackage{pdfcomment}
\usepackage[shortlabels]{enumitem}
\usepackage[english]{babel}
\usepackage{graphicx} 
\usepackage[cmex10]{amsmath}
\usepackage{amsfonts,amssymb,mathtools}
\usepackage{amstext}
\usepackage{upgreek}
\usepackage[tight,footnotesize]{subfigure}
\usepackage[tight,nice]{units}
\usepackage[makeroom]{cancel}
\usepackage[draft]{changes}
\usepackage{import}
\usepackage{tikz}
\usepackage{circuitikz}
\usepackage{pgfplots}
\usepgfplotslibrary{groupplots}
\usepackage{pgfplotstable}
\usepackage{filecontents}
\usepackage{mathtools}
\usepackage{soul}
\usepackage[normalem]{ulem}
\usepackage{color}
\usepackage{blindtext}
\usepackage{tcolorbox}
\usepackage{accents}
\usepackage{bbm}
\usepackage{algorithmicx}
\usepackage{algorithm}
\usepackage{algpseudocode}
\usepackage{collectbox}
\usepackage{caption}
\usepackage{siunitx}

\graphicspath{{./}{../files/fig/}{../files/images/}{../../files/fig/}{../../files/images/}{../../../files/fig/}{../../../files/images/}}

\algnewcommand\algorithmicinput{\textbf{INPUT:}}
\algnewcommand\INPUT{\item[\algorithmicinput]}
\algnewcommand\algorithmicoutput{\textbf{OUTPUT:}}
\algnewcommand\OUTPUT{\item[\algorithmicoutput]}

\numberwithin{equation}{section}

\numberwithin{thm}{section}

\numberwithin{prob}{section}

\graphicspath{{figures/}}

\newcommand{\bs}{\mathbf}

\newcommand{\one}{\bs{1}}

\newcommand{\bx}{\mbox{\boldmath$x$}}

\newcommand{\bn}{\mbox{\boldmath$n$}}
\newcommand{\bu}{\mbox{\boldmath$u$}}
\newcommand{\bv}{\mbox{\boldmath$v$}}

\newcommand{\ba}{\mbox{\boldmath$a$}}
\newcommand{\bb}{\mbox{\boldmath$b$}}

\newcommand{\bd}{\mbox{\boldmath$d$}}
\newcommand{\be}{\mbox{\boldmath$e$}}
\newcommand{\bff}{\mbox{\boldmath$f$}}

\newcommand{\bh}{\mbox{\boldmath$h$}}
\newcommand{\bj}{\mbox{\boldmath$j$}}
\newcommand{\bm}{\mbox{\boldmath$m$}}
\newcommand{\bp}{\mbox{\boldmath$p$}}
\newcommand{\bq}{\mbox{\boldmath$q$}}
\newcommand{\bt}{\mbox{\boldmath$t$}}

\newcommand{\bL}{\mbox{\boldmath$L$}}

\newcommand{\bBB}{\mbox{\boldmath$B$}}

\newcommand{\bDD}{\mbox{\boldmath$D$}}
\newcommand{\bEE}{\mbox{\boldmath$E$}}
\newcommand{\bFF}{\boldsymbol{F}}
\newcommand{\bGG}{\boldsymbol{G}}
\newcommand{\bHH}{\boldsymbol{H}}
\newcommand{\bJJ}{\mbox{\boldmath$J$}}
\newcommand{\bLL}{\boldsymbol{L}}

\newcommand{\bNN}{\boldsymbol{N}}
\newcommand{\bPP}{\boldsymbol{P}}

\newcommand{\bSS}{\boldsymbol{S}}

\newcommand{\bVV}{\boldsymbol{V}}
\newcommand{\bWW}{\boldsymbol{W}}
\newcommand{\bXX}{\boldsymbol{X}}
\newcommand{\bZZ}{\boldsymbol{Z}}
\newcommand{\bAAA}{\boldsymbol{\mathcal{A}}}
\newcommand{\bBBB}{\boldsymbol{\mathcal{B}}}
\newcommand{\bDDD}{\boldsymbol{\mathcal{D}}}
\newcommand{\bEEE}{\boldsymbol{\mathcal{E}}}

\newcommand{\bHHH}{\boldsymbol{\mathcal{H}}}
\newcommand{\bJJJ}{\boldsymbol{\mathcal{J}}}
\newcommand{\bMMM}{\boldsymbol{\mathcal{M}}}
\newcommand{\bPPP}{\boldsymbol{\mathcal{P}}}

\newcommand{\bSSS}{\boldsymbol{\mathcal{S}}}

\newcommand{\bepsilon}{\mbox{\boldmath$\epsilon$} }
\newcommand{\bkappa}{\mbox{\boldmath$\kappa$} }

\newcommand{\bnu}{\mbox{\boldmath$\nu$} }
\newcommand{\bPhi}{\mbox{\boldmath$\Phi$} }
\newcommand{\bsigma}{\mbox{\boldmath$\sigma$} }

\newcommand{\bTheta}{\mbox{\boldmath$\Theta$} }

\newcommand{\bvarphi}{\mbox{\boldmath$\varphi$} }
\newcommand{\bvartheta}{\mbox{\boldmath$\vartheta$} }

\newcommand{\bone}{\boldsymbol{\mathbbm{1}}}

\newcommand{\Hone}[2][]{H^1#1(#2)}
\newcommand{\bHone}[2][]{\boldsymbol{H}^1#1(#2)}

\newcommand{\Hcurlm}[2][]{\boldsymbol{H}#1(\CurlSymbM;#2)}

\newcommand{\GradSymb}{\mathrm{\boldsymbol{grad}}}
\newcommand{\CurlSymb}{\mathrm{\boldsymbol{curl}}}
\newcommand{\DivSymb}{\mathrm{div}}
\newcommand{\GradSymbM}{\text{\bf Grad}}
\newcommand{\CurlSymbM}{\text{\bf Curl}}
\newcommand{\DivSymbM}{\text{Div}}
\newcommand{\Curl}[2][]  {\mathrm{\CurlSymb}{#1}\,{#2}}
\newcommand{\Grad}[2][]  {\mathrm{\GradSymb}{#1}\,{#2}}
\newcommand{\Div}[2][]   {\mathrm{\DivSymb}{#1}\,{#2}}
\newcommand{\Gradm}[2][]  {\text{\GradSymbM}{#1}\,{#2}}
\newcommand{\Curlm}[2][]  {\text{\CurlSymbM}{#1}\,{#2}}
\newcommand{\Divm}[2][]  {\text{\DivSymbM}{#1}\,{#2}}

\journal{Computer Methods in Applied Mechanics and Engineering }

\begin{document}
\begin{frontmatter}
\title{Modeling and simulation of nonlinear electro-thermo-mechanical continua with application to  shape memory polymeric medical devices}
\author{Innocent Niyonzima\fnref{mymainaddress}}
\author{Yang Jiao\fnref{mymainaddress}}
\author{Jacob Fish\fnref{mymainaddress}}
\address[mymainaddress]{Columbia University, Department of Civil Engineering and Engineering Mechanics, New York, 10027 NY, USA.}
\begin{abstract}
	Shape memory materials have gained considerable attention thanks to their ability to change physical properties when subjected to external stimuli such as temperature, pH, humidity, electromagnetic fields, etc. These materials are increasingly used for a large number of biomedical applications. For applications inside the human body, contactless control can be achieved by the addition of electric and/or magnetic particles that can react to electromagnetic fields, thus leading to a composite biomaterial. The difficulty of developing accurate numerical models for smart materials results from their multiscale nature and from the multiphysics coupling of involved phenomena. This coupling involves electromagnetic, thermal and mechanical problems. This paper contributes to the multiphysics modeling of a shape memory polymer material used as a medical stent. The stent is excited by electromagnetic fields produced by a coil which can be wrapped around a failing organ. In this paper we develop large deformation formulations for the coupled electro-thermo-mechanical problem using the electric potential to solve the electric problem. The formulations are then discretized and solved using the finite element method. Results are validated by comparison with results in the literature.
\end{abstract}
\end{frontmatter}
\linenumbers
%
%
\section{Introduction}
\label{sec:motivation}

The increase of life expectancy creates a need to maintain the functions of aging organs to allow greater independence for the elderly. Biomaterial implants have the potential to fulfill some of these functions. The total number of implants in the world exceeds four hundred million per year and grows every year \cite{uweb-biometarials-04}. Biomaterials are also increasingly used for a large number of biomedical applications such as the prevention and cure of coronary heart disease and stroke, as well as ophthalmological applications, biosensors and drug delivery systems \cite{bhatia-biomaterials-10, rezaie-biomaterials-15}. They have the potential to contribute to the reduction of the cost of health and the improvment of the life conditions.

Among biomaterials, shape memory materials have gained considerable attention in the biomedical community thanks to their ability to change physical properties (morphing, structural rigidity, refractive index, etc.) when subjected to external stimuli such as temperature, pH, humidity, electromagnetic fields, etc. This special behavior results from the the shape memory effect observed in shape memory materials \cite{huang-smp-11}. They are used in minimally invasive surgery as embolic devices to treat aneurysm \cite{small-biomaterials-07, maitland-biomaterials-07} and as vascular stents \cite{yakacki-stents-07,baer-stents-09,ajili-stents_09}.
\begin{figure}[!h]
\centering
    \includegraphics[width=1.0\textwidth]{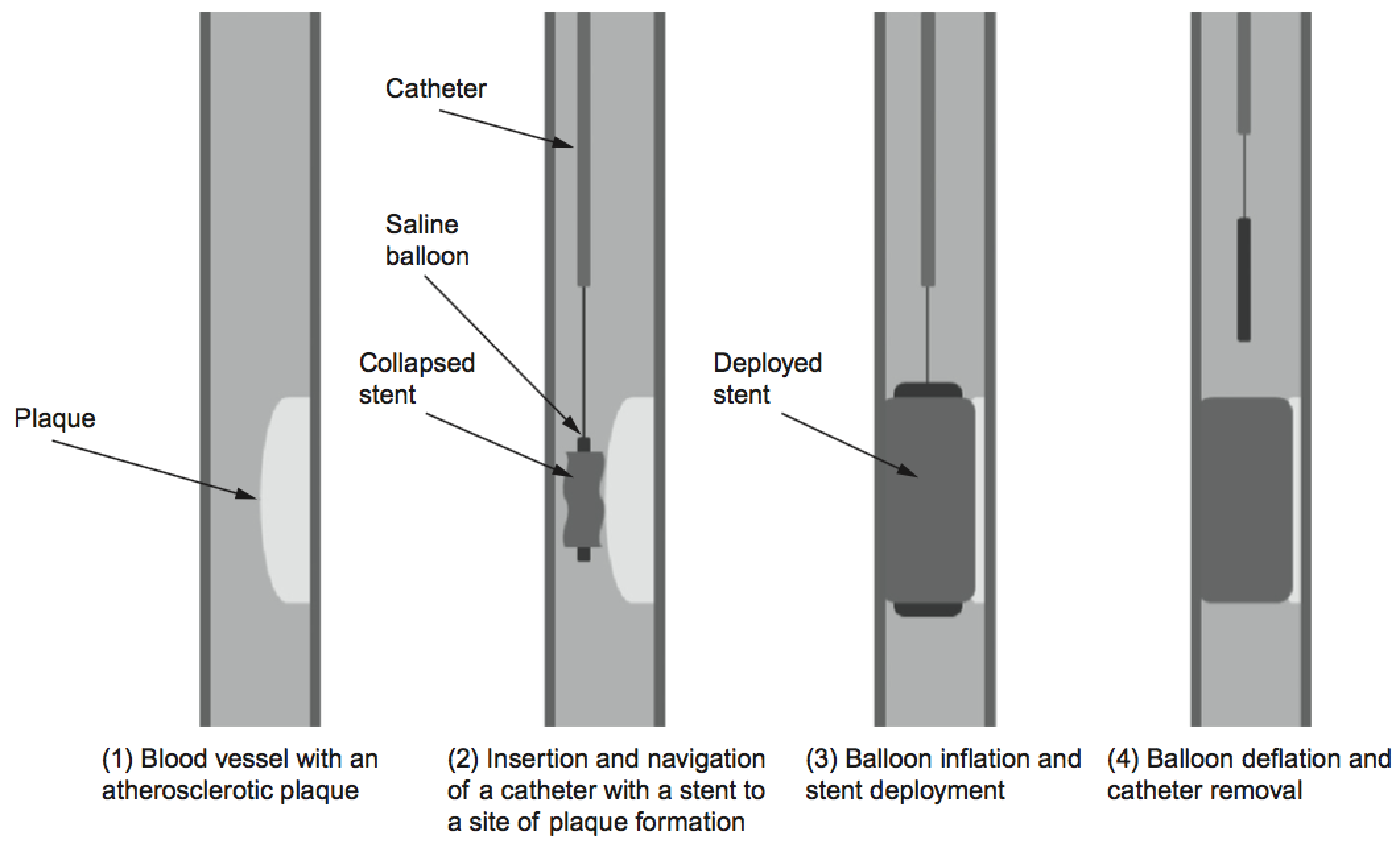}
    \caption{\small Diagram illustrating the deployment of a stent in a blood vessel \cite{yahia-smp-15}.}
\label{fig:sme-stents}
\end{figure}
Figure \ref{fig:sme-stents} illustrate the deployment of a stent in a blood vessel. They can also be used as portable sensors to monitor heart and respiratory rates and in controlled drug delivery systems thus allowing to reduce the side effects of drugs \cite{wischke-drugdelivery-10,nagahama-biomaterial-09}. For these different uses, biomaterials must possess a number of properties. They must be biocompatible to avoid toxicity in contact with biological tissues. Biodegradability is a desirable property for temporary implants, and for minimally invasive in vivo surgery applications, devices must be controllable without contact and self-expanding. All of these properties make polymers the best candidates for a wide range of biomedical applications. Contactless control can be achieved by the addition of electric/magnetic (nano)particles inclusions to produce a smart composite that can react to electromagnetic fields.

Accurate numerical models for smart composites must account for the multiphysics coupling which involve different domains of physics (electromagnetism, thermal, mechanics) and the multiscale nature of the materials. The present manuscript is concerned with multiphysics modeling of shape memory polymer materials used for biomedical devices for homogeneous materials. The difficulty of developing multiphysics models arises from a number of factors. (a) The \emph{geometric non-linearity} resulting from large mechanical deformations leads to the modification of the equations that govern the electromagnetic and thermal problems to account for motion. Additional complexity for the electromagnetic problem results from the presence of electromagnetic fields in the air and vacuum. (b) The \emph{material non-linearities} resulting from the presence of materials with nonlinear thermal constitutive laws, plastic and viscous laws for mechanics and nonlinear anhysteretic/hysteretic constitutive laws for electromagnetism. 

Models and numerical simulations have already been developed for multiphysics problems in piezoelectric, magnetostrictive, and piezomagnetic materials in the case of small deformations \cite{fish-piezo-03,elhadrouz-piezo-06,anderson-coupling-07,khalaquzzaman-piezoelectric-multiscale-12,kuznetsov-hmm-12,perevertov-coupling-15,bishay-piezo-electro-magnetic-15}. Theoretical models have also been developed for thermomechanical and electro-magneto-thermomechanical problems in the case of large mechanical deformations \cite{eringen-electrodynamics-12,pao-electrodynamics-78,ogden-coupling-09,saxena-coupledlargedefo-13}. Popular numerical implementations combine the Lagrangian approach for the mechanical problem and the Eulerian or arbitrary Lagrangian Eulerian approach for the electromagnetic problem \cite{stiemer-ale-09, abali-largedefo-18-a}. In this paper we consider low frequency electromagnetic problems and solve for a scalar potential formulation only defined in the mechanical domain. Thus, a Lagrangian mesh can also be used for the electromagnetic problem. 

The development of multiphysics models for smart composites controlled by electromagnetic fields is still in its infancy. Electro-magneto-mechanical models using electrostatic and magnostatic formulations have been developed for magneto-sensitive composites and magneto-electro-elastic composites (e.g., electro active polymers ) in large deformations \cite{ethiraj-coupling-16,miehe-coupling-16,bayat-coupling-18}. In \cite{homsi-DG-17}, the discontinuous Galerkin method was used to solve the electro-thermo-mechanical problem in a SMP by solving an electrokinetic problem excited by surface currents. For a more effective contactless control, the multiphysics problem should include eddy currents and hysteretic losses as a means of controlling the temperature.

In this paper, we develop a simple multiphysics model for an electromagnetically controlled vascular stent excited by a coil. The paper extends the thermomechanical model developed in \cite{boatti-smp-16} by proposing a contactless electromagnetic control of the temperature, especially during the recovery step that takes place inside the human body. For the sake of clarity and in order to have a self-sufficient paper which is easily accessible by the mechanics and electromagnetic communities, we derive the general fully coupled problem from Maxwell equations and conservation laws using the Lagrangian and Eulerian formalisms. Then, a simplified, quasistatic electro-thermo-mechanical problem is derived and discretized using the finite element (FE) method.

The paper is organized as follows: in Section \ref{section:governing-equations} we recall Maxwell's equations and conservation laws using the Lagrangian and the Eulerian formalisms. In Section \ref{section:formulations}, we derive the simplified coupled problem and its strong and weak forms using potential formulations. The weak formulations are then semi-discretized in space using the FE method and in time using the backward Euler time stepping method. The resulting system of nonlinear algebraic equations is linearized and solved using the Newton--Raphson method. Section \ref{section:numerical_tests} deals with numerical examples. At first, we validate the thermo-mechanical formulation for SMP materials with simple geometries along the lines of \cite{boatti-smp-16}. We then study the behavior of the electromagnetically responsive SMP stent excited by a coil. In Section \ref{section:conclusions} we close the paper with conclusions and perspectives.
%
%
\section{Governing equations of the general multiphysics problem}
\label{section:governing-equations}

In this section, the general electro-magneto-thermo-mechanical coupled problem is derived from Maxwell's equations and conservation laws. Throughout the paper, we use the indices $E$ and $L$ to denote the Eulerian and Lagrangian quantities. Thus, $\bff_E$ and $\bff_L$ denote the forces in the Eulerian and Lagrangian framework, respectively. The open domains $\Omega_{\mathrm{0}}^{\mathrm{Mec}}$, $\Omega_{\mathrm{0}}^{\mathrm{The}}$ and $\Omega_{\mathrm{0}}^{\mathrm{Ele}}$ denote the undeformed computational domains for the mechanical, thermal and electromagnetic problems, respectively. Likewise, $\Omega_{\mathrm{t}}^{\mathrm{Mec}}$, $\Omega_{\mathrm{t}}^{\mathrm{The}}$ and $\Omega_{\mathrm{t}}^{\mathrm{Ele}}$ denote the deformed computational domains for the mechanical, thermal and electromagnetic problems at a time $t \in \mathcal{I}_t := ]t_0, t_{\mathrm{end}}[$. The domains of the mechanical and thermal problems are generally subdomains of the electromagnetic domain, i.e., $\Omega_i^{\mathrm{Mec}} \subseteq \Omega_i^{\mathrm{Ele}}$ and $\Omega_i^{\mathrm{The}} \subseteq \Omega_i^{\mathrm{Ele}}$ with $i = \{0, t\}$ as the electromagnetic fields can be defined in the entire domain including the surrounding air. The domains $\Omega_{c, i} \subseteq \Omega_{i}^{\mathrm{Ele}}$, $\Omega_{c, i}^C \subseteq \Omega_{i}^{\mathrm{Ele}}$ and $\Omega_{s, i} \subsetneq \Omega_{c, i}^C \subset \Omega_{i}^{\mathrm{Ele}}$ are the conductors, non-conductors and inductors where the electric currents source is imposed. The domains $\Gamma_i^{\mathrm{Ele}}$, $\Gamma_i^{\mathrm{Ele}}$ and $\Gamma_i^{\mathrm{Ele}}$ denote the boundaries of the electromagnetic, thermal and mechanical domains, respectively. The differential operators $\Gradm[]{}$, $\Curlm[]{}$ and $\Divm[]{}$ denote the gradient, rotational and divergence operators defined on the undeformed configurations while $\Grad[]{}$, $\Curl[]{}$ and $\Div[]{}$ denote the same operators defined on the deformed configurations.
%
%
\subsection{Kinematics} 

The motion is described by the mappings $\boldsymbol{\bvarphi}_t$ and $\boldsymbol{\bvarphi}$ assumed to be smooth enough (we do not consider fracture). The mapping $\boldsymbol{\bvarphi}_t$ is also assumed to be bijective and defined by: 
\begin{equation}
	\begin{aligned}
		&\bvarphi_t: \,\, && \Omega_0^{\mathrm{Mec}} \rightarrow \mathbb{E}^3, \\
                   & &&\boldsymbol{X} \mapsto \bx = \bvarphi_t(\boldsymbol{X}) = \bvarphi(\boldsymbol{X}, t) = \boldsymbol{X} + \bu(\boldsymbol{X}, t)
	\end{aligned}
	\label{eq:Large_Transformation_Mapping}
\end{equation}
where $\boldsymbol{X}$ is the position of a particle point P in the undeformed configuration, $\bx$ is the position of P in the deformed configuration, $\bu$ is the vector of displacements and $\mathbb{E}^3$ is the three dimensional Euclidean space \cite{wriggers-fem-08}. 
The positions in the undeformed and deformed configurations are related by $\boldsymbol{X} = \bvarphi_t^{-1}(\bx)$ which is valid thanks to the bijection of $\bvarphi_t$. For any time $t \in \mathcal{I}_t$, the deformed configurations are also defined as: 
\begin{equation}
	\Omega_{\mathrm{t}}^{\mathrm{Mec}} := \varphi_t\left(\Omega_{\mathrm{0}}^{\mathrm{Mec}}\right), \quad 
	\Omega_{\mathrm{t}}^{\mathrm{The}} := \varphi_t\left(\Omega_{\mathrm{0}}^{\mathrm{The}}\right),  \quad
	\Omega_{\mathrm{t}}^{\mathrm{Ele}} := \varphi_t\left(\Omega_{\mathrm{0}}^{\mathrm{Ele}}\right).
\end{equation} 

The deformation gradient tensor and its determinant are given by:
\begin{equation}
	\bFF := \frac{\partial \bx}{\partial \bXX} = \Gradm[]{\varphi} = \bone + \Gradm[]{\bu} \quad, 
	\quad J = \mathrm{det}\bFF
    \label{eq:Deformation_Gradient_Tensor}
\end{equation}
where $\bone$ is the identity matrix. The velocity $\bv$ and acceleration $\ba$ are given by:
\begin{gather}
    \bv(\boldsymbol{X}, t) = \frac{\partial \bvarphi}{\partial t}(\boldsymbol{X}, t), \quad 
    \ba(\boldsymbol{X}, t) = \frac{\partial \bv}{\partial t}(\bXX, t) = \frac{\partial^2 \bvarphi}{\partial t^2}(\boldsymbol{X}, t).
    \label{eq:Velocity_Acceleration_Eulerian}
\end{gather}

Assuming the existence of a mapping $\bTheta$:
\begin{equation}
	\begin{aligned}
		&\bTheta: \,\, && \mathbb{E}^3 \times \mathcal{I}_t \rightarrow \Omega_{0}^{\mathrm{Mec}} \subsetneq \mathbb{E}^3, \\
                   & &&(\bx, t) \mapsto \bXX = \bTheta(\bx, t) = \bTheta(\bvarphi(\bXX, t), t),
	\end{aligned}
	\label{eq:Large_Transformation_Inverse_Mapping}
\end{equation}
it is possible to derive the following relationship:
\begin{gather}
    \frac{D \bXX}{D t} 
    = \frac{\partial \boldsymbol{\Theta}}{\partial \bvarphi} \frac{\partial \bvarphi}{\partial t} 
    + \frac{\partial \boldsymbol{\Theta}}{\partial t} 
	= \bFF^{-1} \bv + \bVV = \boldsymbol{0},
\label{eq:Velocity_Lagrangian}
\end{gather}
which relates the \emph{matter flow field} $\bVV$ to the velocity $\bv$ as:
\begin{gather}
    \bVV = -\bFF^{-1} \bv.
\label{eq:Velocity_Lagrangian_V_v}
\end{gather}
The matter flow field is important for the definition of the electromagnetic problem on the undeformed configuration.
The independence of the initial position $\bXX$ on the time was used to derive \eqref{eq:Velocity_Lagrangian}.
%
%
\subsection{Maxwell's equations} 
\label{sub:section_maxwell_equations}

We recall non-relativistic Maxwell's equations on moving domains under large deformations, using Eulerian and Lagrangian formalisms \cite{penfield-electrodynamics-63,fano-electromagnetism-68,pao-electrodynamics-78,eringen-electrodynamics-12}. 

\subsubsection{Full Maxwell's equations in Eulerian formalism} 

In Eulerian setting, the electromagnetic fields are governed by the following Maxwell's equations \cite{jackson-electrodynamics-98,eringen-electrodynamics-12}:
\begin{subequations}
	\begin{equation}
    	\Curl[]{\bh} = \bj + \partial_t \bd, \quad  
    	\Curl[]{\be} = -\partial_t \bb, \quad 
    	\Div[]{\bb} = 0, \quad 
    	\Div[]{\bd} = \rho_{E}, 
	\tag{\theequation\,a-d}
	\end{equation}
	\label{eq:Full_Maxwell_GoverningEquations_Eulerian} 
\end{subequations}
and constitutive laws:
\begin{subequations}
	\begin{multline}
		\bh = \mu_0^{-1} \bb - \underset{\bm_{\mathrm{eff}}}{\underbrace{(\bm - \bv \times \bp)}} = \boldsymbol{\nu}_E(\bb) \, \bb + \bv \times \bp = \bHH(\be, \bb, \bv), \\
		\bd = \epsilon_0 \be + \bp = \epsilon_0 \boldsymbol{\epsilon}_{r E}(\be) \, \be 
		= \boldsymbol{\epsilon}_E(\be) \, \be = \bDD(\be), \textcolor{white}{Inno Inno Inno Inno}\\
        \bj = \bsigma_E (\underset{\be_{\mathrm{eff}}}{\underbrace{\be + \bv \times \bb } } ) + \bj_s + \rho_{E} \bv 
        = \bj_{\mathrm{eff}} + \rho_{E} \bv
        = \bJJ(\be, \bb, \bv).
	\tag{\theequation\,a-c}
	\end{multline}
\label{eq:Full_Maxwell_MaterialLaw_Eulerian}
\end{subequations}
In (\ref{eq:Full_Maxwell_GoverningEquations_Eulerian}\!\,\,a-d)--(\ref{eq:Full_Maxwell_MaterialLaw_Eulerian}\!\,\,a-c), $\bh$ is the magnetic field (A/m), $\bb$ the magnetic flux density (T), $\bj$ the electric current density (A/m$^2$), $\bd$ the electric flux density (C/m$^2$), $\be$ the electric field (V/m) and $\rho_{E}$ the electric charge density (C/m$^3$). The fields $\be_{\mathrm{eff}}$ and $\bj_{\mathrm{eff}}$ are the effective electric field and effective electric current density, whereas $\bj_s$ (A/m$^{2}$) is the electric current source defined in the inductors $\Omega_{s, t}$ and $\bsigma_E$ is the electric conductivity tensor ($\Omega$/m).  
The magnetization $\bm$ (A/m) and polarization $\bp$ (C/m$^2$) are defined by:
\begin{subequations}
	\begin{equation}
		\bm = \mu_0^{-1} \, \boldsymbol{\chi}_{b E}(\bb)\, \bb, \quad 
		\bp = \epsilon_0 \, \boldsymbol{\chi}_{e E}(\be) \, \be 
		= \frac{\boldsymbol{\chi}_{e E}}{\bone + \boldsymbol{\chi}_{e E}} \, \bd, 
		\tag{\theequation\,a-b}
	\end{equation}
\label{eq:Full_Maxwell_MaterialLaw_mb_pe_Eulerian}
\end{subequations}
where $\boldsymbol{\chi}_{b E}$ and $\boldsymbol{\chi}_{e E}$ are the magnetic and electric susceptibility tensors, $\mu_0 = 4 \pi 10^{-7}$ is the magnetic permeability of the free space (H/m) and $\epsilon_0 \simeq 10^{-9}/36 \pi$ is the electric permittivity of the free space (C$^2$/Nm$^2$). Another definition of magnetic susceptibility $\boldsymbol{\chi}_{m E}$ with $\bm = \boldsymbol{\chi}_{m E} \bh$ is often used in the constitutive law dual to (\ref{eq:Full_Maxwell_MaterialLaw_mb_pe_Eulerian}\!\,\,a).  Additionally, $\boldsymbol{\nu}_E = \mu_0^{-1} (\bone - \boldsymbol{\chi}_{b E})$ is the magnetic reluctivity tensor, $\boldsymbol{\mu}_E = \boldsymbol{\nu}_E^{-1}$ is the magnetic permeability tensor, $\bepsilon_E$ is the electric susceptibility tensor, $\bepsilon_{r E}(\be)$ is the relative electric permittivity tensor with $\bepsilon_{r E}(\be) = \bone + \boldsymbol{\chi}_{e E}(\be)$ and $\bv$ is the velocity (m/s). The dependency of the mapping $\bHH$ on (\ref{eq:Full_Maxwell_MaterialLaw_Eulerian}\!\,\,a) on the electric field $\be$ results from (\ref{eq:Full_Maxwell_MaterialLaw_mb_pe_Eulerian}\!\,\,b). The motion is accounted for by the velocity terms in the constitutive laws (\ref{eq:Full_Maxwell_MaterialLaw_Eulerian}\!\,\,a-c).
%
%
\subsubsection{Full Maxwell's equations in Lagrangian formalism} 

In Lagrangian setting, the electromagnetic fields are governed by the following Maxwell's equations (see \cite{penfield-electrodynamics-63,fano-electromagnetism-68,pao-electrodynamics-78,eringen-electrodynamics-12,saxena-coupledlargedefo-13} and \cite[Appendix F]{castro-electromagnetism-14}):
\begin{subequations}
	\begin{equation}
    	\Curlm[]{\bHHH_{\mathrm{eff}}} = \bJJJ + \partial_t \bDDD, \,\, 
    	\Curlm[]{\bEEE_{\mathrm{eff}}} = -\partial_t \bBBB, \,\,
    	\Divm[]{\bBBB} = 0, \,\,
    	\Divm[]{\bDDD} = \rho_{L}, 
	\tag{\theequation\,a-d}
	\end{equation}
	\label{eq:Full_Maxwell_GoverningEquations_Lagrangian} 
\end{subequations}
and constitutive laws:
\begin{align}
	\bHHH_{\mathrm{eff}} &= \underset{\boldsymbol{\nu}_L}{\underbrace{J^{-1} \bFF^{T} (\boldsymbol{\nu}_E \circ \varphi_t^{-1}) \bFF } } \, \bBBB 
	+ (\bepsilon_{r E} \circ \varphi_t^{-1})^{-1} (\bVV \times \bDDD)
	\label{eq:Full_Maxwell_MaterialLaw_Lagrangian_H} \\
	\bDDD &= \underset{\bepsilon_L}{\underbrace{J \, \bFF^{-1} (\bepsilon_E \circ \varphi_t^{-1}) \bFF^{-T}}} (\underset{\bEEE}{\underbrace{\bEEE_{\mathrm{eff}} + \bVV \times \bBBB}}), 
	\label{eq:Full_Maxwell_MaterialLaw_Lagrangian_D}\\
    \bJJJ &= \underset{\boldsymbol{\bsigma}_L}{\underbrace{J \, \bFF^{-1} \, (\bsigma_E \circ \varphi_t^{-1}) \, \bFF^{-T} }} \, \bEEE_{\mathrm{eff}} + \bJJJ_s
	\label{eq:Full_Maxwell_MaterialLaw_Lagrangian_J}
\end{align}
defined on the undeformed configuration $\Omega_{0}^{\mathrm{Ele}}$.
In (\ref{eq:Full_Maxwell_GoverningEquations_Lagrangian}\!\,\,a-d) and \eqref{eq:Full_Maxwell_MaterialLaw_Lagrangian_H}--\eqref{eq:Full_Maxwell_MaterialLaw_Lagrangian_J}, $\bHHH$ is the magnetic field (A/m), $\bBBB$ the magnetic flux density (T), $\bJJJ$ the electric current density (A/m$^2$), $\bDDD$ the electric flux density (C/m$^2$), $\bEEE_{\mathrm{eff}}$ the effective electric field, $\bEEE$ the electric field (V/m), $\bJJJ_s$ (A/m$^{2}$) the current source defined in the inductors $\Omega_{s, 0}$ and $\bVV$ is the matter flow field (m/s) defined in \eqref{eq:Velocity_Lagrangian}. The tensor $\bsigma_L$ is the electric conductivity tensor ($\Omega$/m) and $\rho_{L}$ is the electric charge density (C/m$^3$). The notation $(\bff \circ {\bvarphi_t}^{-1})(\bx) := \bff({\bvarphi_t}^{-1}(\bx))$ is used.

The one differential forms are transformed as: 
\begin{subequations}
	\begin{multline}
    	\bHHH_{\mathrm{eff}} = \bFF^{T} \left(\bh - \bv \times \bd
	 \right), \, \, 
    	\bEEE_{\mathrm{eff}} = \bFF^{T} \left(\be + \bv \times \bb \right), \\
    	\bEEE = \bEEE_{\mathrm{eff}} + \bVV \times \bBBB = \bFF^{T} \be .
	\tag{\theequation\,a-c}
	\end{multline}
\label{eq:Simplified_Maxwell_Pullback_H_E}   
\end{subequations}
and the two differential forms are transformed as:
\begin{equation}
    \bBBB = J \, \bFF^{-1} \, \bb, \quad
    \bDDD = J \, \bFF^{-1} \, \bd, \quad
    \bJJJ = J \, \bFF^{-1} \, \bj, \quad
    \bJJJ_s = J \, \bFF^{-1} \, \bj_s.
\label{eq:Simplified_Maxwell_Pullback_B_D_J}   
\end{equation}
In (\ref{eq:Simplified_Maxwell_Pullback_H_E}\!\,\,c), $\bVV \times \bBBB = -\bFF^{T} \left(\bv \times \bb \right)$ results from the identity:
\begin{equation}
	\bFF^{T} \bv \times \bFF^{T} \bb = J \bFF^{-1}(\bv \times \bb),
\end{equation} 
which is valid for any matrix $\bFF \in GL_3(\mathbb{R})$ \cite[Formula B.11]{castro-electromagnetism-14}.

Additionally, the magnetization $\bMMM$ is related to the magnetic induction by:
\begin{multline}
    \bMMM 
    = \bFF^{T} \bm
    = \bFF^{T} (\mu_0^{-1} \boldsymbol{\chi_{b E}} \, \bb)
    = \bFF^{T} (\mu_0^{-1} \boldsymbol{\chi}_{b E} \, J^{-1} \, \bFF \, \bBBB) \\
    = \mu_0^{-1} J^{-1} \bFF^T \boldsymbol{\chi_{b E}} \bFF \, \bBBB 
    = \mu_0^{-1} \boldsymbol{\chi}_{b L} \, \bBBB
    = ( \mu_0^{-1} \bone - \boldsymbol{\nu}_{L}) \, \bBBB,
\label{eq:Simplified_MaterialLaw_M}   
\end{multline}
while the effective magnetization $\bMMM_{\mathrm{eff}}$ is given by:
\begin{equation}
    \bMMM_{\mathrm{eff}} 
    = \bFF^{T} (\bm - \bv \times \bp)
    = \bMMM + \bVV \times \bPPP,
    = ( \mu_0^{-1} \bone - \boldsymbol{\nu}_{L}) \, \bBBB + \bVV \times \bPPP
\label{eq:Simplified_MterialLaw_M_eff}   
\end{equation}
where the polarization is related to the electric flux density by \cite{castro-electromagnetism-14}:
\begin{equation}
    \bPPP 
    = \frac{\bepsilon_{r L} - \bone}{\bepsilon_{r L}} \, \bDDD
    = \frac{\boldsymbol{\chi_{e E}}}{\bone + \boldsymbol{\chi_{e E}}} \, \bDDD.
\label{eq:Simplified_MterialLaw_M}   
\end{equation}
Combining all these results, the following transformations for second order tensors used in constitutive laws can be derived: 
\begin{multline}
	\boldsymbol{\nu}_{L} = J^{-1} \bFF^{T} \boldsymbol{\nu}_{E} \bFF, 
	\boldsymbol{\mu}_{L} = J \bFF^{-1} \boldsymbol{\mu}_{E} \bFF^{-T}, \\
	\bepsilon_{L} = J \bFF^{-1} \bepsilon_{E} \bFF^{-T}, 	
	\bsigma_{L} = J \bFF^{-1} \bsigma_{E} \bFF^{-T}.
\label{eq:Simplified_Maxwell_Pullback_2ndOrderTensor} 
\end{multline}
%
%
\subsection{Conservation equations} 

We recall conservation equations using the Eulerian and Lagrangian formalisms.

\subsubsection{Conservation equations using the Eulerian description} 

In Eulerian setting, the conservation equations read \cite{simo-fem-06,wriggers-fem-08,belytschko-fem-13}:
\begin{subequations}
    \begin{equation}
        \frac{D \rho}{D t} + \rho \, \Div[]{\bv} = 0, \quad
        \rho \frac{D \bv}{D t} - \Div[]{\bsigma} = \bff_E, \quad
        \bepsilon_{P} \colon \bsigma + \bL_E = 0.
        \tag{\theequation\,a-c}     
    \end{equation}
    \label{eq:Full_Conservation_Equations_Eulerian}     
\end{subequations}
\begin{equation}
    \rho \, \frac{D U_E}{D t} + \Div[]{\bq_E} = \bsigma : \bLL + w_E.
    \label{eq:Full_Conservation_Energy_Eulerian}    
\end{equation}
Equations (\ref{eq:Full_Conservation_Equations_Eulerian}\!\,\,a-c) are balance equations of the mass, linear and angular momentum and \eqref{eq:Full_Conservation_Energy_Eulerian} is the balance equation of the internal energy. The quantity $\rho$ is the mass density (kg\,m$^{-3}$), $\bsigma$ is the Cauchy stress (N/m$^2$), $\bff_E$ is the volume force (N/m$^3$), $\bepsilon_{P}$ is the third order permutation tensor also known as the Levi Civita tensor such that ${(\bepsilon_p \bsigma)}_i = (\bepsilon_p)_{ijk} \bsigma_{jk}$, $\bL_E$ is the torque, $U_E$ is the density of internal energy, $\bq_E$ is the heat flux density, $\bLL := \Grad[]{\bv}$ is the gradient deformation and $w_E$ is the electromagnetic source term for the heat problem. 

The electromagnetic force, torque, internal energy and source term are given by: 
\begin{multline}
    \textcolor{white}{I\,\,} \bff_E = \rho_{E} \be + \bj \times \bb + \left(\Grad[]{\be} \right)^{T} \bp + 
    \left(\Grad[]{\bb} \right)^{T} \bm \\+ \partial_t (\bp \times \bb) + \Div[]{[\bv \otimes (\bp \times \bb)]},
    \label{eq:Full_Electromagnetic_Force_Eulerian}
\end{multline}
\begin{align}
    \bL_E &= \bp \times \be + (\bm - \bv \times \bp) \times \bb, 
    	\textcolor{white}{\bL_E = \bp \times \be + (\bm - \bv \times \bp) \times \bb. }
		\label{eq:Full_Electromagnetic_Torque_Eulerian} \\
    \rho \, \frac{D U_E}{D t} &= \rho \, c_p \frac{\partial \vartheta_E}{\partial t}, 
    	\textcolor{white}{\bL_E = \bp \times \be + (\bm + \bv \times \bp) \times \bb, } 
			\label{eq:Full_Electromagnetic_Energy_Eulerian} \\
    w_E &= \bj_{\mathrm{eff}} \cdot \be_{\mathrm{eff}}  
		- \bm_{\mathrm{eff}} \cdot \frac{\partial \bb}{\partial t}
    	+ \rho \frac{\partial}{\partial t} \left( \frac{\bp}{\rho} \right) \cdot \be_{\mathrm{eff}} 
    	\label{eq:Full_Electromagnetic_Losses_Eulerian}
\end{align}
where $c_p$ and $\vartheta_E$ are the heat capacity and the temperature, respectively. Equations (\ref{eq:Full_Conservation_Equations_Eulerian} c) and \eqref{eq:Full_Electromagnetic_Torque_Eulerian} suggest that Cauchy stress $\bsigma$ may not be symmetric in presence of electromagnetic fields. However, symmetry may be kept for isotropic materials.

Equations \eqref{eq:Full_Conservation_Equations_Eulerian}--\eqref{eq:Full_Conservation_Energy_Eulerian} must be completed by constitutive laws derived from the \emph{Clausius-Duhem inequality}. In this paper, we assume the following nonlinear constitutive law for the thermal problem \cite{belytschko-fem-13}:
\begin{equation}
	\bq_E = \boldsymbol{Q}(\vartheta_E, \Grad[]{\vartheta_E}) = -\bkappa_E(\vartheta_E) \, \Grad[]{\vartheta_E}
\label{eq:Full_Thermal_Mechanical_ConstitutiveLaws_Eulerian}
\end{equation}
where $\bkappa_E$ is the thermal conductivity tensor (W/mK). The thermo-mechanical constitutive law involves the definition of an appropriate objective rate $\bsigma^{\nabla}$ (e.g., Jaumann rate, Truesdell rate or Green-Naghdi rate) which is related to the material derivative of the Cauchy stress $\dot{\bsigma}$ \cite{belytschko-fem-13} as:
\begin{subequations}
	\begin{equation}
		\bsigma^{\nabla} = \dot{\bsigma} - \boldsymbol{\psi}(\boldsymbol{L}, \bFF) 
		= \boldsymbol{\Sigma}(\bsigma, \boldsymbol{L}, \bFF, \vartheta_E, \bZZ_E(\tau \leq t)), \quad 
		\boldsymbol{L} = \dot{\bFF} \, \bFF^{-1}
		\tag{\theequation\,a-b}
	\end{equation}
\label{eq:Full_ThermoMechanical_ConstitutiveLaws_Eulerian}
\end{subequations}
where the velocity gradient $\bLL$ is related to the rate of the deformation gradient $\dot{\bFF}$ through (\ref{eq:Full_ThermoMechanical_ConstitutiveLaws_Eulerian}\!\,\,b) and $\bZZ_E(\tau \leq t)$ is the set of internal variables that account for the history of the loadings.
%
%
\subsubsection{Conservation equations using the Lagrangian description} 

In Lagrangian setting, the conservation equations read \cite{simo-fem-06,wriggers-fem-08,belytschko-fem-13}:
\begin{subequations}
    \begin{equation}
        \rho_0 - \rho J = 0, \quad
        \rho_0 \frac{D^2 \bu}{D t^2} - \Divm[]{\bFF \bSS} = \bff_L, \quad
        \bepsilon_{P} \colon (\bFF \bSS \bFF^{T}) + \bL_L = 0.
        \tag{\theequation\,a-c}     
    \end{equation}
    \label{eq:Full_Conservation_Equations_Lagrangian}     
\end{subequations}
\begin{equation}
    \rho_0 \, \frac{D U_L}{D t} + \Div[]{\bq_L} = -\bPP^{T} \colon \Gradm[]{(\bFF \bVV)} + w_L.
    \label{eq:Full_Conservation_Energy_Lagrangian}    
\end{equation}
Equations (\ref{eq:Full_Conservation_Equations_Lagrangian}\!\,\,a-c) are balance equations of the mass, linear and angular momentum and \eqref{eq:Full_Conservation_Energy_Lagrangian} is the balance equation of the internal energy, expressed on the undeformed configuration $\Omega_0^{\mathrm{Mec}}$. 
The quantity $\rho_0$ is the mass density (kg/m$^{3}$), $\bSS$ is the second Piola--Kirchhoff stress (N/m$^2$) with $\bSS = \bFF^{-1} \bPP$ where $\bPP$ is the first Piola--Kirchhoff stress or nominal stress tensor with $\bPP = J \, \bsigma \, \bFF^{-T}$. The quantity $\bff_L$ is the volume force (N/m$^3$), $\bL_L$ is the torque, $U_L$ is the density of internal energy, $\bq_L$ is the heat flux density and $w_L$ is the source term for the heat equation. The electromagnetic force which is a three differential form is transformed as:
\begin{equation}
	\begin{aligned}
		\!\!\bff_L \!=\! J \, \bff_E 
		&= \rho_L \bFF^{-T} \bEEE + J^{-1} \left( (\bFF \, \bJJJ) \times (\bFF \, \bBBB) \right)  \\
    	&+ J \left(\bFF^{-T} \Gradm[]{(\bFF^{-T} \, \bEEE)} \right)^{T} (J^{-1} \, \bFF \, \bPPP) \\
    	&+ \left(\bFF^{-T} \Gradm[]{(J^{-1} \, \bFF \, \bBBB)} \right)^{T} (\bFF^{-T} \, \bMMM) \\
    	&+ J \frac{\partial}{\partial t}\left[ J^{-2} (\bFF \bPPP) \times (\bFF \bBBB) \right] 
    	+\Divm[]{\left[J^{-1} \bVV \otimes \left[ (\bFF \bPPP) \times (\bFF \bBBB) \right] \right]}.
	\end{aligned}
	\label{eq:Full_Electromagnetic_Force_Lagrangian}
\end{equation}
The torque and the internal energy are transformed as:
\begin{subequations}
	\begin{align}
	    \!\!\!\!\!\!\!\!\!\!\!\!\bL_L 
	    &= J \, \bL_E = \left(\bFF \bPPP \right) \times \left(\bFF^{-T} \bEEE \right) 
		+ \left(\bFF^{-T} \bMMM_{\mathrm{eff}} \right) \times \left(\bFF \bBBB \right), \\
		\!\!\!\!\!\!\!\!\!\!\!\!\rho_0 \frac{D U_L}{D t} 
		&= \rho_0 c_p \frac{\partial \vartheta_L}{\partial t}, 
		\textcolor{white}{\bL_E = \bp \times \be + (\bm + \bv \times \bp) \times \bb, } 
	\end{align}
\label{eq:Full_Electromagnetic_Torque_Losses_Lagrangian}
\end{subequations}
\noindent where $\vartheta_L$ is the temperature expressed on the undeformed configuration and the source term is obtained using the transformation: 
\begin{subequations}
	\begin{align}
	w_L = J \, w_E 
		&= (\bFF \bJJJ) \cdot (\bFF^{-T} \bEEE_{\mathrm{eff}} ) \\
	    &-J \left(\bFF^{-T} \bMMM_{\mathrm{eff}} \right) \times \left( \partial_t (J^{-1} \bFF \bBBB) + 	
	 	\Gradm[]{ \left( J^{-1} \bFF \bBBB \right) } \bVV \right) \\
		&- \rho_0 \left( \frac{\partial}{\partial t} \left(\frac{\bFF \bPPP}{\rho_0} \right) 
		+ \Gradm[]{\left( \frac{\bFF \bPPP}{\rho_0}\right) \bVV}\right) \cdot \left( \bFF^{-T} \bEEE_{\mathrm{eff}} \right).
	\end{align}
\label{eq:Full_Electromagnetic_Losses_Lagrangian}
\end{subequations}

Equations \eqref{eq:Full_Conservation_Equations_Lagrangian}--\eqref{eq:Full_Conservation_Energy_Lagrangian} must be completed by constitutive laws which relate the stress tensor to its associated conjugate strain tensor. In this paper, we use the second Piola--Kirchhoff stress and the Green--Lagrange strain tensors. We assume the following nonlinear constitutive laws for the thermal and thermo-mechanical problems \cite{belytschko-fem-13}:
\begin{multline}
	\bq_L 
		= J \, \bFF^{-1} \, \boldsymbol{Q}(\Gradm[]{(\vartheta_E \circ {\bvarphi_t}^{-1})}) \\
		= \underset{\bkappa_L}{\underbrace{J \, \bFF^{-1} \, \bkappa_E \,  \bFF^{-T} } } \Gradm[]{\vartheta_L}
		= -\bkappa_L(\vartheta_L) \, \Gradm[]{\vartheta_L},
\label{eq:Full_Thermal_ConstitutiveLaws_Lagrangian}
\end{multline}
\begin{subequations}
	\begin{equation}
		\bSS = \bSSS_{\mathrm{VEP}}(\bEE, \dot{\bEE}, \vartheta_L, \bZZ_L(\tau \leq t)), \quad 
		\bEE = \frac{1}{2} (\bFF^T \bFF - \bone),
		\tag{\theequation\,a-b}
	\end{equation}
\label{eq:Full_ThermoMechanical_ConstitutiveLaws_Lagrangian}
\end{subequations}
where $\bkappa_L$ is the thermal conductivity tensor (W/mK) and $\bSSS_{\mathrm{VEP}}$ is a mapping that represents the visco-elastoplastic constitutive law. Viscosity is reflected through the dependence of the stress on the rate of the Green--Lagrange tensor $\dot{\bEE}$, plasticity is accounted for using a set of internal variables $\bZZ_L(\tau \leq t)$ and the thermo-mechanical aspect is accounted for by the dependency on the temperature $\vartheta_L$.
%
%
%
%
\section{Formulations, discretization and linearization}
\label{section:formulations}

In this section, we derive the simplified multiphysics problem from the fully coupled problem defined in section \ref{section:governing-equations}. Using this simplified problem, strong and weak formulations of the multiphysics problems are derived using potentials. The weak formulations are then discretized in space using the continuous Galerkin approximation and in time using the backward Euler integrator. Finally, the resulting nonlinear system of algebraic equations is linearized. 

\subsection{Simplified governing equations of the multiphysics problem}
\label{section:simplified_governing-equations}

We make the following \emph{magnetoquasistatic (MQS) assumptions}

\begin{equation}
	\delta_i \simeq L_{\mathrm{sys, i}}, \quad \lambda \gg L_{\mathrm{sys, i}}.
    \label{eq:Simplified_Maxwell_GoverningEquations_MQSAssumption}
\end{equation}
In \eqref{eq:Simplified_Maxwell_GoverningEquations_MQSAssumption}, $\delta_i := 1/\sqrt{\pi f \sigma_{E, i} \mu_{E, i}}$ is the skin depth in the spatial direction $i = x, y$ and $z$, $f$ is the frequency of the source term, $\sigma_{E, i}$ and $\mu_{E, i}$ are eigenvalues of $\bsigma_E$ and $\boldsymbol{\mu}_E$,  $\lambda$ is the wavelength corresponding to the frequency $f$ and $L_{\mathrm{sys, i}}$ is the characteristic length of the structure in a spatial direction $i$ \cite{hiptmair-mqs-05}. The first term of \eqref{eq:Simplified_Maxwell_GoverningEquations_MQSAssumption} explains the presence of eddy currents while the second explains the neglect of electromagnetic waves. Further in this section, we will relax the first condition in \eqref{eq:Simplified_Maxwell_GoverningEquations_MQSAssumption} to $\delta_i = \alpha \, L_{\mathrm{sys, i}}$ with $\alpha$ which is big thus allowing to neglect the reaction field. 

\subsubsection{The magnetoquasistatic problem} 

Using the MQS assumption, the following MQS problem in the deformed configuration can be defined: 
\begin{subequations}
    \begin{equation}
        \Curl[]{\bh} = \bj, \quad  
        \Curl[]{\be} = -\partial_t \bb, \quad 
        \Div[]{\bb} = 0, 
        \tag{\theequation\,a-c}
    \end{equation}
    \label{eq:Simplified_Maxwell_GoverningEquations_Eulerian}
\end{subequations}
together with the constitutive laws:
\begin{subequations}
	\begin{equation}
        \bh = \boldsymbol{\nu}_E(\bb) \, \bb = \bHH(\bb), \quad
        \bj = \bj_{\mathrm{eff}} = \bsigma_E \be_{\mathrm{eff}} + \bj_s = \bsigma_E (\be + \bv \times \bb) + \bj_s.
    \tag{\theequation\,a-b}
    \end{equation}
    \label{eq:Simplified_Maxwell_MaterialLaw_Eulerian} 
\end{subequations}

In Lagrangian setting, the MQS problem is governed by Maxwell's equations \cite{eringen-electrodynamics-12,castro-electromagnetism-14}:
\begin{subequations}
    \begin{equation}
        \Curlm[]{\bHHH_{\mathrm{eff}}} = \bJJJ, \quad
        \Curlm[]{\bEEE_{\mathrm{eff}}} = -\partial_t \bBBB, \quad
        \Divm[]{\bBBB} = 0
    \tag{\theequation\,a-c} 
    \end{equation}
\label{eq:Simplified_Maxwell_GoverningEquations_Lagrangian}    
\end{subequations}
completed by the following constitutive laws:
\begin{subequations}
	\begin{equation}
    	\bHHH_{\mathrm{eff}} = \nu_{L} \, \bBBB, \quad
    	\bJJJ = \bsigma_{L} \, \bEEE_{\mathrm{eff}} + \bJJJ_s.
    \tag{\theequation\,a-b} 
	\end{equation}
\label{eq:Simplified_Maxwell_MaterialLaw_Lagrangian_H_J}   
\end{subequations}

\subsubsection{Conservation equations} 

In addition to the MQS assumption, we assume \emph{quasistatic mechanical problem} thus neglecting the inertia term in the balance of linear momentum and \emph{isotropic magnetic materials} therefore restoring the symmetry of the Cauchy and the second Piola--Kirchhoff stress as $\bLL_E = \bLL_L = \boldsymbol{0}$ in (\ref{eq:Full_Conservation_Equations_Eulerian} c) and (\ref{eq:Full_Conservation_Equations_Lagrangian} c). Additionally, we neglect mechanical losses in the heat equation. Under these assumptions, balance equations in the deformed configuration become:
\begin{subequations}
    \begin{equation}
        \frac{D \rho}{D t} + \rho \Div[]{\bv} = 0, \quad
        \Div[]{\bsigma} + \bff_E = 0, \quad
        \bsigma = \bsigma^T.
        \tag{\theequation\,a-c}     
    \end{equation}
    \label{eq:Simplified_Conservation_Equations_Eulerian}   
\end{subequations}
\begin{equation}
    \rho \, c_p \frac{\partial \vartheta_E}{\partial t} + \Div[]{\bq_E} = w_E 
    \label{eq:Simplified_Conservation_Energy_Eulerian}   
\end{equation}
and the electromagnetic force, torque and electromagnetic losses in \eqref{eq:Full_Electromagnetic_Force_Eulerian}--\eqref{eq:Full_Electromagnetic_Losses_Eulerian} become:
\begin{subequations}
	\begin{equation}
    	\bff_E = \bj \times \bb + \left(\Grad[]{\bb} \right)^{T} \bm, \quad 
    	\bL_E = \bm  \times \bb, \quad
    	w_E = \bj_{\mathrm{eff}} \cdot \be_{\mathrm{eff}} - \bm_{\mathrm{eff}} \cdot \frac{\partial \bb}{\partial t}.
    	\tag{\theequation\,a-c}     
	\end{equation}
	\label{eq:Simplified_Electromagnetic_Torque_Losses_Eulerian} 
\end{subequations}

In Lagrangian setting, the conservation equations become:
\begin{subequations}
    \begin{equation}
        \rho_0 - \rho J = 0, \quad
        \Divm[]{\bFF \bSS} + \bff_L = \boldsymbol{0}, \quad
        \bSS = \bSS^T,
        \tag{\theequation\,a-c}
    \end{equation}
    \label{eq:Simplified_Conservation_Equations_Lagrangian}   
\end{subequations}
\begin{equation}
    \rho_0 \, c_p \, \frac{\partial \vartheta_L}{\partial t} + \Divm[]{\bq_L} = w_L.
    \label{eq:Simplified_Conservation_Energy_Lagrangian}   
\end{equation}
where the electromagnetic force and torque are given by:
\begin{align}
    \bff_L &= J^{-1} \left( (\bFF \, \bJJJ) \times (\bFF \, \bBBB) \right) + 
    J\left(\bFF^{-T} \Gradm[]{(J^{-1} \, \bFF \, \bBBB}) \right)^{T} (\bFF^{-T} \, \bMMM), 
	\label{eq:Simplified_Electromagnetic_Force_Lagrangian} \\
    \bL_L &= \left(\bFF^{-T} \bMMM_{\mathrm{eff}} \right) 
    \times \left(\bFF \bBBB \right),
	\label{eq:Simplified_Electromagnetic_Torque_Lagrangian} 
\end{align}
and the source term is given by:
\begin{multline}
    	w_L = (\bFF \bJJJ ) \cdot (\bFF^{-T} \bEEE_{\mathrm{eff}} ) - \\
    	J (\bFF^{-T} \bMMM_{\mathrm{eff}}) \times \left[ \partial_t (J^{-1} \bFF \bBBB) + \Gradm[]{ \left( J^{-1} \bFF \bBBB \right) } \bVV \right].
\label{eq:Simplified_Electromagnetic_Losses_Lagrangian}
\end{multline}
Further in the paper, we consider elasto-plastic materials governed by the following constitutive law:
\begin{subequations}
	\begin{equation}
		\bSS = \bSSS_{EP}(\bEE, \vartheta_L, \bZZ_L(\tau \leq t)), \quad 
		\bEE = \frac{1}{2} (\bFF^T \bFF - \bone).
		\tag{\theequation\,a-b}
	\end{equation}
\label{eq:Simplified_ThermoMechanical_ConstitutiveLaws_Lagrangian}
\end{subequations}
In this case, the second Piola--Kirchhoff stress does not depend on the rate of Green--Lagrange strain as viscosity is not considered.

\subsection{Strong forms} 

\subsubsection{The magnetoquasistatic problem} 

The strong formulation of the MQS problem is derived using the so-called \emph{magnetic induction conforming formulations} \cite{bossavit-cem-98}. In the deformed configuration, the derivation is achieved by verifying equations (\ref{eq:Simplified_Maxwell_GoverningEquations_Eulerian}\!\,\,b-c) in the strong sense:
\begin{equation}
    \bb = \Curl[]{\ba} \simeq \bb_s = \Curl[]{\ba_s}, \, \,
    \be = -\partial_t \ba - \Grad[]{\phi} \simeq -\partial_t \ba_s - \Grad[]{\phi}
    \label{eq:Strong_Derivation_Electromagnetic_Potentials_Eulerian}
\end{equation}
where $\ba$ is the vector potential unknown field used for the eddy currents problem, $\ba_s$ is the source vector potential which can be pre-computed based on the electric current source $\bj_s$ imposed in inductors such as in coils and $\phi$ is the unknown scalar potential. The approximation $\ba \simeq \ba_s$ implies the neglect of the reaction field  and is valid under the assumption $\delta_i = \alpha \, L_{\mathrm{sys, i}}$ with $\alpha$ which is big \cite{scorretti-formulations-12}. 

The electric potential $\phi$ is governed by the following problem where \eqref{eq:Strong_Electrokinetic_Governing_Equation_Eulerian} is obtained by applying the divergence $\Div[]{}$ operator to (\ref{eq:Simplified_Maxwell_GoverningEquations_Eulerian} a) and \eqref{eq:Strong_Electrokinetic_Constitutive_Law_Eulerian} is derived from (\ref{eq:Simplified_Maxwell_MaterialLaw_Eulerian} b):
\begin{align}
    \Div[]{\bj} &= 0 && \text{ in } \Omega_{\mathrm{t}}^{\mathrm{Ele}}, 
    \label{eq:Strong_Electrokinetic_Governing_Equation_Eulerian} \\
    \bj &= \bsigma_{E} (\be + \bv \times \bb) + \bj_s &&\\
        &= \bsigma_{E} (-\partial_t \ba_s - \Grad[]{\phi} + \bv \times \Curl[]{\ba_s}) + \bj_s
    	&& \text{ in } \Omega_{\mathrm{t}}^{\mathrm{Ele}},
    	\label{eq:Strong_Electrokinetic_Constitutive_Law_Eulerian} \\
    \phi(\bx, t) &= \phi_{D}(\bx, t) 
    && \text{ on } \Gamma_t^{\mathrm{Diri, Ele}},
    \label{eq:Strong_Electrokinetic_Dirichlet_BC_Eulerian} \\
    \bn_E \cdot \bj &= 0
    && \text{ on } \Gamma_t^{\mathrm{Neu, Ele}}.
    \label{eq:Strong_Electrokinetic_Neumann_BC_Eulerian}
\end{align}
\eqref{eq:Strong_Electrokinetic_Dirichlet_BC_Eulerian} and \eqref{eq:Strong_Electrokinetic_Neumann_BC_Eulerian} are the Dirichlet and the Neumann boundary conditions.

In the undeformed configuration, the derivation is carried out by verifying equations (\ref{eq:Simplified_Maxwell_GoverningEquations_Lagrangian}\!\,\,b-c) in the strong sense:
\begin{multline}
    \bBBB = \Curlm[]{\bAAA} \simeq \bBBB_s = \Curlm[]{\bAAA_s} = \Curlm[]{\left(\bFF^{T} \ba_s\right)}, \\
    \bEEE_{\mathrm{eff}} = -\partial_t \bAAA - \Gradm[]{\Phi} \simeq -\partial_t \bAAA_s - \Gradm[]{\Phi} 
     = -\bFF^{T} \partial_t \ba_s - \Gradm[]{\Phi}.
     \label{eq:Strong_Derivation_Electromagnetic_Potentials_Eulerian}
\end{multline}
In \eqref{eq:Strong_Derivation_Electromagnetic_Potentials_Eulerian}, the 1 differential form $\bAAA_s$ is transformed as $\bAAA_s = \bFF^{T} \ba_s$. The electric potential $\Phi$ is therefore governed by the following problem where \eqref{eq:Strong_Electrokinetic_Governing_Equation_Lagrangian} is obtained by applying the divergence operator $\Divm[]{}$ to (\ref{eq:Simplified_Maxwell_GoverningEquations_Lagrangian}\!\,\,a) and \eqref{eq:Strong_Electrokinetic_Constitutive_Law_Lagrangian} is derived from (\ref{eq:Simplified_Maxwell_MaterialLaw_Lagrangian_H_J}\!\,\,b):
\begin{align}
    \Divm[]{\bJJJ} &= 0, && \text{ in } \Omega_0^{\mathrm{Ele}}
    \label{eq:Strong_Electrokinetic_Governing_Equation_Lagrangian}, \\
    \bJJJ 
    &= \bsigma_{L} \, \bEEE_{\mathrm{eff}} + \bJJJ_s \\
    &= -\bsigma_{L} (\bFF^{T} \partial_t \ba_s + \Gradm[]{\Phi}) + \bJJJ_s , && \text{ in } \Omega_0^{\mathrm{Ele}}
    \label{eq:Strong_Electrokinetic_Constitutive_Law_Lagrangian} \\
    \Phi(\bx, t) &= \Phi_{D}(\bx, t) 
    && \text{ on } \Gamma_0^{\mathrm{Diri, Ele}}
    \label{eq:Strong_Electrokinetic_Dirichlet_BC_Lagrangian }, \\
    \bn_L \cdot \bJJJ &= 0 
    && \text{ on } \Gamma_0^{\mathrm{Neu, Ele}}
    \label{eq:Strong_Electrokinetic_Neumann_BC_Lagrangian}.
\end{align}
where the conductivity tensor is transformed as $\bsigma_{L} = J \bFF^{-1} \bsigma_{E}(\vartheta_E \circ \varphi_t^{-1}) \bFF^{-T}$ thanks to \eqref{eq:Simplified_Maxwell_Pullback_2ndOrderTensor}.

\subsubsection{The heat equation} 

In the deformed configuration, the evolution of the temperature is governed by the following problem derived from \eqref{eq:Simplified_Conservation_Energy_Eulerian} and (\ref{eq:Simplified_Electromagnetic_Torque_Losses_Eulerian} c) and \eqref{eq:Full_Thermal_ConstitutiveLaws_Lagrangian}:
\begin{align}
    \rho c_{p} \frac{\partial \vartheta_E}{\partial t} + \Div[]{\bq_E} &= w_E 
    &&\text{ in } \Omega_{\mathrm{t}}^{\mathrm{The}},
    \label{eq:Strong_Heat_Governing_Equation_Eulerian} \\
    \bq_E &= -\bkappa_E(\vartheta_E) \, \Grad[]{\vartheta_E}      
    && \text{ in } \Omega_{\mathrm{t}}^{\mathrm{The}},
    \label{eq:Strong_Heat_Constitutive_Law_Eulerian} \\
    \vartheta_E(\bx, 0) &= \vartheta_{E,0}(\bx) 
    && \text{ in } \Omega_0^{\mathrm{The}}
    \label{eq:Strong_Heat_IC_Eulerian}, \\
    \vartheta_E(\bx, t) &= \vartheta_{E, D}(\bx, t) 
    && \text{ on } \Gamma_t^{\mathrm{Diri, The}}
    \label{eq:Strong_Heat_Dirichlet_BC_Eulerian}, \\
    \bn_E \cdot \bq_E &= h_E(t)(\vartheta_E - \vartheta_{E, B}) 
    && \text{ on } \Gamma_t^{\mathrm{conv, The}}
    \label{eq:Strong_Heat_Convective_BC_Eulerian}, \\
    \bn_E \cdot \bq_E &= \epsilon_E^R \sigma_E^R (\vartheta_E^4 - \vartheta_{E, R}^4) 
    && \text{ on } \Gamma_t^{\mathrm{rad, The}}
    \label{eq:Strong_Heat_Radiative_BC_Eulerian}.
\end{align}
The source term in \eqref{eq:Strong_Heat_Governing_Equation_Eulerian} can be expressed in terms of the potentials as:
\begin{multline}
	w_E = \underset{w_E^{\mathrm{eddy}}}{\underbrace{\left(\bsigma_E \left(\partial_t \ba_s + \Grad[]{\phi} - \bv \times \Curl[]{\ba_s} \right) \right)^2}} \\
		\underset{w_E^{\mathrm{eddy}}}{\underbrace{- \bj_s \cdot (\bsigma_E \, (\partial_t \ba_s + \Grad[]{\phi} - \bv \times \Curl[]{\ba_s}))}}
		- \underset{w_E^{\mathrm{hyst}}}{\underbrace{\bm_{\mathrm{eff}} \cdot \Curl[]{ \left(\frac{\partial \ba_s}{\partial t} \right) } } }
\end{multline}
where $w_E^{\mathrm{eddy}}$ and $w_E^{\mathrm{hyst}}$ represent eddy current and hysteretic losses.
Equations \eqref{eq:Strong_Heat_IC_Eulerian}, \eqref{eq:Strong_Heat_Dirichlet_BC_Eulerian}, \eqref{eq:Strong_Heat_Convective_BC_Eulerian} and \eqref{eq:Strong_Heat_Radiative_BC_Eulerian} represent the initial condition, Dirichlet, convective and radiative boundary conditions, respectively.

The Lagrangian strong form of the heat problem is given by: 
\begin{align}
    \rho_0 c_{p} \frac{\partial \vartheta_L}{\partial t} + \Divm[]{\bq_L} &= w_L 
    &&\text{ in } \Omega_{\mathrm{0}}^{\mathrm{The}}
    \label{eq:Strong_Heat_Governing_Equation_Lagrangian}, \\
    \bq_L &= -\bkappa_L(\vartheta_L) \, \Gradm[]{\vartheta_L}     
    && \text{ in } \Omega_{\mathrm{0}}^{\mathrm{The}}
    \label{eq:Strong_Heat_Constitutive_Law_Eulerian}, \\
    \vartheta_L(\bx, 0) &= \vartheta_{L,0}(\bx) 
    && \text{ in } \Omega_0^{\mathrm{The}}
    \label{eq:Strong_Heat_IC_Lagrangian}, \\
    \vartheta_L(\bx, t) &= \vartheta_{L, D}(\bx, t) 
    && \text{ on } \Gamma_0^{\mathrm{Diri, The}}
    \label{eq:Strong_Heat_Dirichlet_BC_Lagrangian}, \\
    \bn_L \cdot \bq_L &= h_L(t)(\vartheta_L - \vartheta_{L, B}) 
    && \text{ on } \Gamma_0^{\mathrm{conv, The}}
    \label{eq:Strong_Heat_Convective_BC_Lagrangian},\\
    \bn_L \cdot \bq_L &= \epsilon_L^R \sigma_L^R (\vartheta_L^4 - \vartheta_{L, R}^4) 
    && \text{ on } \Gamma_t^{\mathrm{rad, The}}
    \label{eq:Strong_Heat_Radiative_BC_Lagrangian}
\end{align}
where the thermal conductivity is transformed according to \eqref{eq:Simplified_Maxwell_Pullback_2ndOrderTensor} as $\bkappa_L = J\, \bFF^{-1} \, \bkappa_E \bFF^{-T}$ and the convective heat coefficient $h_L(t)$ is transformed using Nanson's formula as $h_L(t) = J \, | \bFF^{-T} \bn_L| \, h_E(t)$.

The source term in \eqref{eq:Strong_Heat_Governing_Equation_Lagrangian} can be expressed in terms of the potentials as:
\begin{multline}
	w_L = (\bFF^{\textcolor{white}{{-T}}} \!\!\!\!\!\!\!\! (\underset{-\bJJJ}{\underbrace{\bsigma_L \left(\partial_t \bAAA_s 
	+ \Gradm[]{\Phi} \right) - \bJJJ_s } } ) )
	\cdot
	(\underset{-\bEEE_{\mathrm{eff}}}{\underbrace{\bFF^{-T} \left(\partial_t \bAAA_s + \Gradm[]{\Phi} \right) } } ) + 
	\\
	J \bFF^{-T} \underset{\bMMM_{\mathrm{eff}}}{\underbrace{\mu_0^{-1}(\bone - \bnu_L) \Curlm[]{\bAAA_s} } } \bigg[ \partial_t (J^{-1} \bFF \Curlm[]{\bAAA_s} ) + \Gradm[]{\left(J^{-1} \bFF \Curlm[]{\bAAA_s} \right)} \bigg].
\end{multline}
Equations \eqref{eq:Strong_Heat_IC_Lagrangian}, \eqref{eq:Strong_Heat_Dirichlet_BC_Lagrangian}, \eqref{eq:Strong_Heat_Convective_BC_Lagrangian} and \eqref{eq:Strong_Heat_Radiative_BC_Lagrangian} represent the initial condition, Dirichlet, convective and radiative boundary conditions, respectively.

\subsubsection{The mechanical problem} 

Mechanical fields in the undeformed configuration are governed by the following problem:
\begin{align}
    \Divm[]{(\bFF \, \bSS)} + \bff_L &= 0 &&\text{ in } \Omega_{\mathrm{0}}^{\mathrm{Mec}},
    \label{eq:Strong_Mechanics_Governing_Equation_Lagrangian} \\
    \bSS &= \bSSS_{EP}\left(\bEE, \vartheta_L, \bZZ_{L}(\tau \leq t) \right) 
    && \text{ in } \Omega_{0}^{\mathrm{Mec}},
    \label{eq:Strong_Mechanics_Constitutive_Law_Lagrangian}\\
    \bEE &= \frac{\bFF^T \bFF - \bone}{2} && \text{ in } \Omega_{0}^{\mathrm{Mec}},
    \label{eq:Strong_GreenLagrange_Deformation}\\
    \bu(\bx, t) &=  \bu_D(\bx, t) && \text{ on } \Gamma_{0, D}^{\mathrm{Mec}}.
    \label{eq:Strong_Mechanics_Dirichlet_BC_Lagrangian}\\
    \bn_L \cdot (\bFF \bSS) &=  \bt_L && \text{ on } \Gamma_{0, N}^{\mathrm{Mec}}.
    \label{eq:Strong_Mechanics_Neumann_BC_Lagrangian}
\end{align}
In terms of the potential, the force $\bff_L$ in \eqref{eq:Simplified_Electromagnetic_Force_Lagrangian} is given by:
\begin{multline}
    \bff_L = \bFF_L(\bu, \Phi) 
    = J^{-1} \, \bFF \, \underset{\bJJJ}{\underbrace{(\bsigma_L (\partial_t \bAAA_s + \Gradm[]{\Phi}) + \bJJJ_s)}} 
    \times
    \bFF \, \underset{\bBBB}{\underbrace{\Curlm[]{\left(\bAAA_s\right)}}} + 
    \\
    J (\bFF^{-T}  \, \Gradm[]{( J^{-1} \, \bFF \, \underset{\bBBB}{\underbrace{\Curlm[]{\left(\bAAA_s\right)}}} ) } )^T 
    (\bFF^{-T} \underset{\bMMM}{\underbrace{\mu_0^{-1} \boldsymbol{\chi_{B L} } J^{-1} \, \bFF^T \, \bFF \, \Curlm[]{\bAAA_s} } } ).
\label{eq:Strong_Mechanics_Force_Lagrangian}
\end{multline}
The term $\bt_L$ in \eqref{eq:Strong_Mechanics_Neumann_BC_Lagrangian} represents the surface traction applied on part of the boundary $\Gamma_{0, N}^{\mathrm{Mec}}$. The thermo-mechanical constitutive law \eqref{eq:Strong_Mechanics_Constitutive_Law_Lagrangian}--\eqref{eq:Strong_GreenLagrange_Deformation} is derived from (\ref{eq:Simplified_ThermoMechanical_ConstitutiveLaws_Lagrangian}\!\,\,a). In this paper, we use the constitutive law described in \cite{boatti-smp-16}.
As a reminder, the total deformation gradient $\bFF$ and the total deformation gradients in the glassy and rubbery states were given by:
\begin{equation}
    \bFF = \bFF^{tg} = \bFF^{tr} = \bone + \Gradm[]{\bu},
    \label{eq:Strong_Mechanics_Constitutive_Law_Lagrangian_Kinematics_1_new}
\end{equation}
where the superscript $t$ denotes the total deformation gradient and the superscripts $r$ and $g$ were used for the rubbery and the glassy states. The total deformation gradients of both phases are decomposed as:
\begin{equation}
    \bFF^{tg} = \bFF^{g}\bFF^{f} = \bFF^{eg}\bFF^{pg}\bFF^{f}, \bFF^{tr} = \bFF^{r}\bFF^{p} = \bFF^{er}\bFF^{p},
    \label{eq:Strong_Mechanics_Constitutive_Law_Lagrangian_Kinematics_2_new}
\end{equation}
where $\bFF^{eg}$ and $\bFF^{pg}$ are the deformation gradients for the elastic and plastic phases in the glassy state, $\bFF^{f}$ is the frozen deformation gradient that represents the temporary deformation which is stored during high temperature shape fixing and $\bFF^{r}$ and $\bFF^{p}$ are the elastic and plastic deformation gradients for the rubbery state.

The total Cauchy stress was also given by:
\begin{equation}
    \bsigma = z^g \bsigma^g + (1 - z^g) \bsigma^r,
    \label{eq:Strong_Mechanics_Constitutive_Law_sigma_new}
\end{equation}
where the temperature-dependent parameter $z^g$ is the ratio of the glassy state. Using the second Piola--Kirchhoff stress and Green-Lagrange strain tensors, the following expression of the stress tensor can be derived:
\begin{equation}
    \begin{aligned}
        \bSS &= J \bFF^{-1} \bsigma \bFF^{-T} = z^g J \bFF^{-1} \bsigma^g \bFF^{-T} + (1 - z^g) J \bFF^{-1} \bsigma^r \bFF^{-T}\\
             &= z^g J^f {\bFF^{f}}^{-1} \underset{\bSS^{g}}{\underbrace{\left(J^g {\bFF^{g}}^{-1} \bsigma^g {\bFF^{g}}^{-T} \right)}} {\bFF^{f}}^{-T} \\
             & \textcolor{white}{Inno Inno Inno}+ (1 - z^g) J^p {\bFF^{p}}^{-1} \underset{\bSS^{r}}{\underbrace{\left(J^r {\bFF^{r}}^{-1} \bsigma^r {\bFF^{r}}^{-T} \right)}} {\bFF^{p}}^{-T}. 
    \end{aligned}
    \label{eq:Strong_Mechanics_Constitutive_Law_sigma_PK_new}
\end{equation}
In \eqref{eq:Strong_Mechanics_Constitutive_Law_sigma_PK_new}, $\bSS^{g}$ and $\bSS^{r}$ are the second Piola--Kirchhoff stress tensors defined on the intermediate configurations \cite{boatti-smp-16} by:
\begin{align}
    &\bSS^g =  \left(\bFF^{pg}\right)^{-1} \left(\lambda^g \mathrm{tr}(\bEE^{eg}) \bone + 2 \mu^g \bEE^{eg} \right) 
    \left(\bFF^{pg}\right)^{-T},
    \label{eq:Strong_Mechanics_Constitutive_Law_Lagrangian_1_g_new} \\
    &\bSS^r =  \lambda^r \mathrm{tr}(\bEE^{er}) \bone + 2 \mu^r \bEE^{er}. 
    \label{eq:Strong_Mechanics_Constitutive_Law_Lagrangian_1_r_new}
\end{align}
The parameters $\lambda^g$, $\lambda^r$, $\mu^g$ and $\mu^r$ are Lam\'{e} parameters for the glassy and the rubbery states. Determinants of deformation gradients of intermediate configurations are defined as $J^{i} = \mathrm{det}\bFF^{i}$ where the superscript $^\emph{i}$ refer to configurations of the glassy state ($^\emph{g}$, $^\emph{eg}$ and $^\emph{pg}$) or of the rubbery state ($^\emph{er}$ and $^\emph{p}$). The contribution to the stress due to thermal expansion have been neglected in \eqref{eq:Strong_Mechanics_Constitutive_Law_Lagrangian_1_g_new}--\eqref{eq:Strong_Mechanics_Constitutive_Law_Lagrangian_1_r_new}.

Details on the equations that govern the evolution of internal variables $\bZZ_{L}(\tau \leq t) := (z^g, \bFF^{f}, \bFF^{p},\bFF^{pg})$ and the numerical update of the internal variables can be found in \cite{boatti-smp-16}.

\subsection{Weak forms} 

The weak forms of the electromagnetic problem \eqref{eq:Strong_Electrokinetic_Governing_Equation_Lagrangian}--\eqref{eq:Strong_Electrokinetic_Neumann_BC_Lagrangian}, the heat problem \eqref{eq:Strong_Heat_Governing_Equation_Lagrangian}--\eqref{eq:Strong_Heat_Radiative_BC_Lagrangian} and the mechanical problem \eqref{eq:Strong_Mechanics_Governing_Equation_Lagrangian}--\eqref{eq:Strong_Mechanics_Neumann_BC_Lagrangian} read \cite{bossavit-cem-98,bachinger-cem-05,wriggers-fem-08,hughes-fem-12}: \\
for each $t \in \mathcal{I}_t$, find $\left(\Phi \times \vartheta_L \times \bu \right) \in U \times V \times \bWW$  such that
\begin{equation}
    \int_{\Omega_0^{\mathrm{Ele}}} \bsigma_L \Gradm[]{\Phi} \cdot \Gradm[]{\Phi^{'}} \, \mathrm{d} \Omega_0 + 
    \int_{\Omega_0^{\mathrm{Ele}}} \bsigma_L \partial_t \bAAA_s \cdot \Gradm[]{\Phi^{'}} \, \mathrm{d} \Omega_0 = 0,
    \label{eq:Weak_Electrokinetic_Equation_Lagrangian}
\end{equation}
\begin{multline}
    \int_{\Omega_0^{\mathrm{The}}} \rho_0 c_p \frac{\partial \vartheta_L}{\partial t}  \cdot \vartheta_{L}^{'} \, \mathrm{d} \Omega_0 + 
    \int_{\Omega_0^{\mathrm{The}}} \underset{-\bq_L}{\underbrace{\bkappa_L \, \Gradm[]{\vartheta_{L}}}}  \cdot \Gradm[]{\vartheta_{L}^{'}} \, \mathrm{d} \Omega_0 
    - \textcolor{white}{Inno Inno Inno Inno Inno Inno Inno Inno}
    \\
    \int_{\Omega_0^{\mathrm{The}}} \underset{w_L}{\underbrace{ (\bFF\left(\bsigma_L \left(\partial_t \bAAA_s + \Gradm[]{\Phi} \right) - \bJJJ_s \right) )
	\cdot
	(\bFF^{-T} \left(\partial_t \bAAA_s + \Gradm[]{\Phi} \right) ) } } \cdot \vartheta_{L}^{'} \, \mathrm{d} \Omega_0
    \\
    +
    \int_{\Gamma_0^{\mathrm{conv, The}}} \underset{\bn_L \cdot \bq_L}{\underbrace{h_L(t)(\vartheta_{L} - \vartheta_{L, B})}} \cdot \vartheta_L^{'} \, \mathrm{d} \Gamma_0
    +
    \\ 
    \int_{\Gamma_0^{\mathrm{rad, The}}} \underset{\bn_L \cdot \bq_L}{\underbrace{ \epsilon_L^R \sigma_L^R (\vartheta_{L}^4 - \vartheta_{L, R}^4)}} \cdot \vartheta_L^{'} \, \mathrm{d} \Gamma_0
     = 0
    \label{eq:Weak_Thermal_Equation_Lagrangian}
\end{multline}
\begin{multline}
    \int_{\Omega_0^{\mathrm{Mec}}} \bSSS_{EP}(\bu, \vartheta_L, \bZZ_{L}) \colon \delta \bEE \, \mathrm{d} \Omega_0 - 
    \\
    \int_{\Omega_0^{\mathrm{Mec}}} \bFF_L(\bu, \Phi) \cdot \bu^{'} \, \mathrm{d} \Omega_0
    =
    \int_{\Gamma_0^{\mathrm{N}}} \bt_L \cdot \bu^{'} \, \mathrm{d} \Gamma_0
\label{eq:Weak_Mechanics_Equation_Lagrangian}
\end{multline}
holds for all test functions $\left(\Phi^{'} \times \vartheta_{L}^{'} \times \bu^{'}\right) \in \left(U_0 \times V_0 \times \bWW_0 \right)$. The force $\bff_L = \bFF_L(\bu, \Phi)$ is given by \eqref{eq:Strong_Mechanics_Force_Lagrangian} and the dependence on the displacement is achieved through $J$ and $\bFF$.
The function spaces are defined such that 
$U \subseteq \Hone[]{\Omega_0^{\mathrm{Ele}}}$, 
$V \subseteq \Hone[]{\Omega_0^{\mathrm{The}}}$, 
$\bWW \subseteq \bHone[]{\Omega^{\mathrm{Mec}}} \equiv \left(\Hone[]{\Omega^{\mathrm{Mec}}}\right)^3$ and the source term of the electromagnetic problem $\bAAA_s$ belongs to a subspace of $\Hcurlm[]{\Omega_0^{\mathrm{Ele}}}$. The virtual Green--Lagrange strain $\delta \bEE$ is related to the virtual displacement $\bu^{'}$ through $\delta \bEE = \bFF^{T} \Gradm[]{\bu^{'}}$.

\subsection{Spatial and temporal discretization}

The unknown fields $\Phi$, $\vartheta_L$ and $\bu$ in \eqref{eq:Weak_Electrokinetic_Equation_Lagrangian}--\eqref{eq:Weak_Mechanics_Equation_Lagrangian} belong to infinite dimensional functional spaces. For numerical simulation, these fields need to be approximated by finite dimensional spaces
\begin{equation}
	\Phi(\bx, t) \approx \bar{\Phi}(\bx, t), \quad 
	\vartheta_L(\bx, t) \approx \bar{\vartheta}_L(\bx, t), \quad
	\bu(\bx, t) \approx \bar{\bu}(\bx, t)
\label{eq:Spatial_Discretization_Lagrangian_Approximation_Spaces}
\end{equation}
defined by:
\begin{equation}
	\bar{\Phi}(\bx, t) = \sum_{i = 1}^{N^{\mathrm{Ele} } } \bar{\Phi}_i(t) N_i^{\mathrm{Ele} }(\bx), \, \, 
	\Gradm[]{\bar{\Phi}}(\bx) = \sum_{i = 1}^{N^{\mathrm{Ele} } } \bar{\Phi}_i(t) \Gradm[]{N_i}^{\mathrm{Ele} }(\bx), 
	\label{eq:Spatial_Discretization_Lagrangian_Ele_Phi}
\end{equation}
\begin{equation}
	\bar{\vartheta}_L(\bx, t) = \sum_{i = 1}^{N^{\mathrm{The} } } \bar{\vartheta}_i(t) N_i^{\mathrm{The} }(\bx), \, 
	\Gradm[]{\bar{\vartheta}_L}(\bx) = \sum_{i = 1}^{N^{\mathrm{The}} } \bar{\vartheta}_i(t) \Gradm[]{N_i^{\mathrm{The} }}(\bx), 
	\label{eq:Spatial_Discretization_Lagrangian_The_VarTheta}
\end{equation}
\begin{equation}
	\bar{\bu}(\bx, t) = \sum_{i = 1}^{N^{\mathrm{Mec} } } \bar{\bu}_i(t) N_i^{\mathrm{Mec} }(\bx), \, \, 
	\Gradm[]{\bar{\bu}}(\bx) = \sum_{i = 1}^{N^{\mathrm{Mec}} } \bar{\bu}_i(t) \Gradm[]{N_i^{\mathrm{Mec} }}(\bx), 
	\label{eq:Spatial_Discretization_Lagrangian_Mec_U} 
\end{equation}
where $N^{\mathrm{Ele}}$, $N^{\mathrm{The}}$ and $N^{\mathrm{Mec}}$ are the number of nodes of the electromagnetic, thermal and mechanical domains, $\bar{\Phi}_i$, $\bar{\vartheta}_i$ and $\bar{\bu}_i = (\bar{u}_{i, x}, \bar{u}_{i, y}, \bar{u}_{i, z})$ are degrees of freedom and $N_i^{\mathrm{Ele} }$, $N_i^{\mathrm{The} }$ and $N_i^{\mathrm{Mec} }$ are shape functions for the electromagnetic, thermal and mechanical problems, respectively.

Inserting \eqref{eq:Spatial_Discretization_Lagrangian_Ele_Phi}--\eqref{eq:Spatial_Discretization_Lagrangian_Mec_U} into \eqref{eq:Weak_Electrokinetic_Equation_Lagrangian}--\eqref{eq:Weak_Mechanics_Equation_Lagrangian} leads to the following discrete system of equations:
\begin{align}
    \displaystyle \boldsymbol{K}^{\mathrm{Ele}}(\bar{\bvartheta}, \bar{\bu}) \bar{\bPhi} 
	+ \bFF^{\mathrm{Ele}}(\bar{\bvartheta}, \bar{\bu}) = \boldsymbol{0}, 
	\label{eq:Spatial_Discretized_Lagrangian_a} \\
    \displaystyle \boldsymbol{M}^{\mathrm{The}} \frac{D \bar{\bvartheta}}{D t} 
	+ \boldsymbol{K}^{\mathrm{The}}(\bar{\bvartheta}, \bar{\bu}) \bar{\bvartheta} 
	+ \bFF^{\mathrm{The}}(\bar{\bPhi}, \bar{\bvartheta}, \bar{\bu}) = \boldsymbol{0}, 
	\label{eq:Spatial_Discretized_Lagrangian_b} \\
    \displaystyle \boldsymbol{K}^{\mathrm{Mec}}(\bar{\bvartheta}, \bar{\bu}, \bZZ_{L}) 
	+ \bar{\bFF}^{\mathrm{Mec}}(\bar{\bPhi}, \bar{\bu}) 
	= \boldsymbol{0}.
	\label{eq:Spatial_Discretized_Lagrangian_c}
\end{align}
where $\bar{\bPhi}$, $\bar{\bvartheta}_L$ and $\bar{\bu}$ are vectors of degrees of freedom and the matrices in \eqref{eq:Spatial_Discretized_Lagrangian_a}--\eqref{eq:Spatial_Discretized_Lagrangian_c} are given by:
\begin{align}
	\boldsymbol{K}^{\mathrm{Ele} } &= \sum_{e = 1}^{N^{\mathrm{Ele} } } {\bLL^{e}}^T 
	\left[ \int_{\Omega_e} (\bBB^{e})^T \, \bsigma_L(\bar{\bu}, \bar{\bvartheta}) \, \bBB^{e} \mathrm{d} \Omega_e \right] {\bLL^{e}}^T,
	\label{eq:Spatial_Lagrangian_Terms_K_Ele} \\
	\bFF^{\mathrm{Ele} } &= \sum_{e = 1}^{N^{\mathrm{Ele} } } {\bLL^{e}}^T 
	\left[ \int_{\Omega_e} (\bBB^{e})^T \, \bsigma_L(\bar{\bu}, \bar{\bvartheta}) \, \bAAA_s^{e} \mathrm{d} \Omega_e \right],
	\label{eq:Spatial_Lagrangian_Terms_F_Ele} 
	\\
	\boldsymbol{M}^{\mathrm{The} } &= \sum_{e = 1}^{N^{\mathrm{The} } } {\bLL^{e}}^T 
	\left[ \int_{\Omega_e} (\bNN^{e})^T \, \rho_0 c_p \, \bNN^{e} \mathrm{d} \Omega_e \right] {\bLL^{e}}^T,
	\label{eq:Spatial_Lagrangian_Terms_M_The} \\ 
	\boldsymbol{K}^{\mathrm{The}} &= \sum_{e = 1}^{N^{\mathrm{The} } } {\bLL^{e}}^T 
	\left[ \int_{\Omega_e} (\bBB^{e})^T \, \bkappa_L(\bar{\bu}, \bar{\bvartheta}) \, \bBB^{e} \mathrm{d} \Omega_e \right] {\bLL^{e}}^T,
	\label{eq:Spatial_Lagrangian_Terms_K_The} 
	\\
	\bFF^{\mathrm{The} } &= \sum_{e = 1}^{N^{\mathrm{The} } } {\bLL^{e}}^T 
	\left[ \int_{\Omega_e} (\bBB^{e})^T \, \bar{w}_L(\bar{\bu}, \bar{\bvartheta}, \bar{\bPhi}) \, \mathrm{d} \Omega_e \right],
	\label{eq:Spatial_Lagrangian_Terms_F_The} \\
	\boldsymbol{K}^{\mathrm{Mec}} &= \sum_{e = 1}^{N^{\mathrm{Mec} } } {\bLL^{e}}^T 
	\left[ \int_{\Omega_e} (\bBB^{e})^T \, \bSSS_{EP}(\bar{\bu}, \bar{\bvartheta}, \bZZ_{L}) \, \mathrm{d} \Omega_e \right], 
	\label{eq:Spatial_Lagrangian_Terms_K_Mec} \\
	\bFF^{\mathrm{Mec} } &= \sum_{e = 1}^{N^{\mathrm{Mec}}} {\bLL^{e}}^T 
	\left[ \int_{\Omega_e} (\bNN^{e})^T \, \bFF_L(\bar{\bu}, \bar{\bPhi}) \mathrm{d} \Omega_e \right]
	\label{eq:Spatial_Lagrangian_Terms_F_Mec} 
\end{align}
where $\bLL^{e}$ is the gather matrix, $\bNN^{e}$ is the element shape function matrix, $\bBB^{e}$ is composed of the elements of the  gradient of $\bNN^{e}$ \cite{fish-fem-07, wriggers-fem-08} and where we neglected the boundary terms in \eqref{eq:Weak_Electrokinetic_Equation_Lagrangian}--\eqref{eq:Weak_Mechanics_Equation_Lagrangian}. The recurrent dependence on the displacement field in \eqref{eq:Spatial_Discretized_Lagrangian_a}--\eqref{eq:Spatial_Discretized_Lagrangian_c} and \eqref{eq:Spatial_Lagrangian_Terms_K_Ele}--\eqref{eq:Spatial_Lagrangian_Terms_F_Mec} results from the transformations between the deformed and undeformed configurations which involve the deformation gradient $\bFF$ and its determinant $J$, and the electric conductivity $\bsigma_L(\bar{\bu}, \bar{\bvartheta})$ is considered to be temperature-dependent.

Equations \eqref{eq:Spatial_Discretized_Lagrangian_a}--\eqref{eq:Spatial_Discretized_Lagrangian_c} can be written as a system of differential algebraic equations:
\begin{equation}
 		\underset{\displaystyle \boldsymbol{M}}{\underbrace{
		\begin{pmatrix}
  			\boldsymbol{0} & \boldsymbol{0}                & \boldsymbol{0} \\
  			\boldsymbol{0} & \boldsymbol{M}^{\mathrm{The}} & \boldsymbol{0} \\
  			\boldsymbol{0} & \boldsymbol{0}                & \boldsymbol{0} 
 		\end{pmatrix}
		}}
		\displaystyle \frac{D }{D t}
		\underset{\displaystyle \bar{\bv}}{\underbrace{
 		\begin{pmatrix}
  			\bar{\boldsymbol{\bPhi}} \\
  			\bar{\boldsymbol{\bvartheta}} \\
  			\bar{\bu}  
 		\end{pmatrix}
		}}	
		\!+\!
		\underset{\displaystyle \bff(\bar{\bv}, \bZZ_{L})}{\underbrace{
 		\displaystyle \begin{pmatrix}
  			\boldsymbol{K}^{\mathrm{Ele}}(\bar{\bvartheta}, \bar{\bu}) \bar{\bPhi} + \bFF^{\mathrm{Ele}}(\bar{\bvartheta}, \bar{\bu}) 
			\\
  			\displaystyle \boldsymbol{K}^{\mathrm{The}}(\bar{\bvartheta}, \bar{\bu}) \bar{\bvartheta} + \bFF^{\mathrm{The}}(\bar{\bPhi}, \bar{\bvartheta}, \bar{\bu}) 
			\\
  			\displaystyle \boldsymbol{K}^{\mathrm{Mec}}(\bar{\bvartheta}, \bar{\bu}, \bZZ_{L}) + \bar{\bFF}^{\mathrm{Mec}}(\bar{\bPhi}, \bar{\bu}) 
 		\end{pmatrix}
		}}
		\!=\! 
 		\begin{pmatrix}
  			\boldsymbol{0} \\
  			\boldsymbol{0} \\
  			\boldsymbol{0} 
 		\end{pmatrix}	 
\label{eq:Spatial_Discretized_Lagrangian_b}
\end{equation}
or
\begin{equation}
 	\displaystyle \boldsymbol{M} \frac{D \bar{\bv} }{D t} + \bff(\bar{\bv}, \bZZ_{L}) = \boldsymbol{0},
\label{eq:DAE_Continuous_Lagrangian}
\end{equation}
where $\boldsymbol{M}$ is a singular matrix and 
\begin{multline}
	\bff_1 = \boldsymbol{K}^{\mathrm{Ele}}(\bar{\bvartheta}, \bar{\bu}) \bar{\bPhi} 
	+ \bFF^{\mathrm{Ele}}(\bar{\bvartheta}, \bar{\bu}),
	\bff_2 = \displaystyle \boldsymbol{K}^{\mathrm{The}}(\bar{\bvartheta}, \bar{\bu}) \bar{\bvartheta} 
	+ \bFF^{\mathrm{The}}(\bar{\bPhi}, \bar{\bvartheta}, \bar{\bu}), \\
	\bff_3 = \displaystyle \boldsymbol{K}^{\mathrm{Mec}}(\bar{\bvartheta}, \bar{\bu}, \bZZ_{L}) 
	+ \bar{\bFF}^{\mathrm{Mec}}(\bar{\bPhi}, \bar{\bu}).
\label{eq:Spatial_Discretized_Lagrangian_Terms}
\end{multline}
Equation \eqref{eq:DAE_Continuous_Lagrangian} can be discretized in time using the backward Euler integrator:
\begin{equation}
 	\displaystyle \boldsymbol{M} \frac{\bar{\bv}^{n+1} - \bar{\bv}^{n}}{\Delta t} + \bff(\bar{\bv}^{n+1}, \bZZ_{L}) 
	= \boldsymbol{0},
\label{eq:DAE_Discretized_Lagrangian}
\end{equation}
where $\bar{\bv}^{n+1} = \bar{\bv}(t_{n+1})$ with $t_{n+1} = t_0 + (n+1) \Delta t$ and $\Delta t$ which is the time step. 
After reordering the terms of \eqref{eq:DAE_Discretized_Lagrangian}, the following nonlinear equations can be derived:
\begin{multline}
 	\displaystyle \boldsymbol{M} \bar{\bv}^{n+1} + {\Delta t} \bff(\bar{\bv}^{n+1}, \bZZ_{L}) - \boldsymbol{M} \bar{\bv}^{n} 
	= \bGG(\bar{\bv}^{n+1}, \bZZ_{L}) - \boldsymbol{M} \bar{\bv}^{n} \\
	= \bHH(\boldsymbol{\bar{\Phi}}^{n+1}, \boldsymbol{\bar{\vartheta}}^{n+1}, \bar{\bu}^{n+1}, \bZZ_{L}) - \boldsymbol{M} \bar{\bv}^{n} 
	= \boldsymbol{0}
\label{eq:Nonlinear_Lagrangian_1}
\end{multline}
with the vector function $\bHH = (\bHH_1, \bHH_2, \bHH_3)$ defined such that $\bHH_i(\boldsymbol{\bar{\Phi}}^{n+1}, \boldsymbol{\bar{\vartheta}}^{n+1}, \bar{\bu}^{n+1}, \bZZ_{L}):= \bGG_i(\bar{\bv}^{n+1}, \bZZ_{L})$ for i = 1, 2 or 3.

\subsection{Linearization}

Equation \eqref{eq:Nonlinear_Lagrangian_1} can be solved using the Newton--Raphson method. To do this, an iterative schema is used with the following linearization:
\begin{multline}
 	\bGG(\bar{\bv}^{n+1}, \bZZ_{L}) \simeq \displaystyle \bGG(\bar{\bv}_{m}^{n+1}, \bZZ_{L})
	+ \displaystyle \left(\frac{\partial \bGG}{\partial \bar{\bv}^{n+1}} \right)_{\bar{\bv}_{m}^{n+1} } 
	\left(\bar{\bv}_{m+1}^{n+1} - \bar{\bv}_{m}^{n+1} \right) 
	\\ 
	= \bb_{\mathrm{RHS}}({\bar{\bv}_{m}^{n+1}}) - \boldsymbol{A}({\bar{\bv}_{m}^{n+1}}) \, \Delta \bar{\bv}_{m+1}^{n+1} = \boldsymbol{0}, 
    \label{eq:Linearization_Lagrangian}
\end{multline}
where the index $m$ is used to denote the Newton--Raphson iteration. The terms in \eqref{eq:Linearization_Lagrangian} are given by:
\begin{multline}
 	\bb_{\mathrm{RHS}} = \bGG(\bar{\bv}_{m}^{n+1}, \bZZ_{L})
	= \bHH(\boldsymbol{\bar{\Phi}}_{m}^{n+1}, \boldsymbol{\bar{\vartheta}}_m^{n+1}, \bar{\bu}_{m}^{n+1}, \bZZ_{L})
	\\
	=
	\begin{pmatrix}
  		\bHH_1(\bar{\bPhi}_{m}^{n+1}, \bar{\bvartheta}_{m}^{n+1}, \bar{\bu}_{m}^{n+1})
		\\
		\bHH_2(\bar{\bPhi}_{m}^{n+1}, \bar{\bvartheta}_{m}^{n+1}, \bar{\bu}_{m}^{n+1})
		\\
		\bHH_3(\bar{\bPhi}_{m}^{n+1}, \bar{\bvartheta}_{m}^{n+1}, \bar{\bu}_{m}^{n+1}, \bZZ_{L})
	\end{pmatrix}.
\label{eq:Linearization_Lagrangian_NL_Term}
\end{multline}
The term $\bHH_3$ depends on the internal variables $\bZZ_L(\tau \leq t)$ through the second Piola--Kirchhoff stress $\bSS$. For each quadrature point, the stress is updated using the return mapping described in \cite{boatti-smp-16}.

The stiffness matrix is given by:
\begin{equation}
 	\!\!\left(\frac{\partial \bGG}{\partial \bar{\bv}^{n+1}} \right)_{\bar{\bv}_{m}^{n+1} }
	\!\!= \boldsymbol{A} = 
	\begin{pmatrix}
		\displaystyle \frac{\partial \bHH_1}{\partial \bar{\bPhi}^{n+1} } &
		\displaystyle \frac{\partial \bHH_1}{\partial \bar{\bvartheta}^{n+1} } &
		\displaystyle \frac{\partial \bHH_1}{\partial \bar{\bu}^{n+1} } 
		\\
		\displaystyle \frac{\partial \bHH_2}{\partial \bar{\bPhi}^{n+1} } &
		\displaystyle \frac{\partial \bHH_2}{\partial \bar{\bvartheta}^{n+1} } &
		\displaystyle \frac{\partial \bHH_2}{\partial \bar{\bu}^{n+1} } 
		\\
		\displaystyle \frac{\partial \bHH_3}{\partial \bar{\bPhi}^{n+1} } &
		\displaystyle \frac{\partial \bHH_3}{\partial \bar{\bvartheta}^{n+1} } &
		\displaystyle \frac{\partial \bHH_3}{\partial \bar{\bu}^{n+1} } 
 	\end{pmatrix}.
    \label{eq:Exact_Linearization_Lagrangian_NL_Term}
\end{equation}
The terms of the tangent stiffness matrix in \eqref{eq:Exact_Linearization_Lagrangian_NL_Term} are given by:
\begin{align}
	\displaystyle \left( \frac{\partial \bHH_1}{\partial \bar{\bPhi}^{n+1} } \right)_{\bar{\bv}_{m}^{n+1} } 
	&= \left(\boldsymbol{K}^{\mathrm{Ele}} \right)_{\bar{\bPhi}_{m}^{n+1}, \bar{\bvartheta}_{m}^{n+1}, \bar{\bu}_{m}^{n+1}},
	\label{eq:Inexact_Linearization_Lagrangian_NL_Term_a_1} \\
	\displaystyle \left( \frac{\partial \bHH_1}{\partial \bar{\bvartheta}^{n+1} } \right)_{\bar{\bv}_{m}^{n+1} } 
	&= \left( \frac{\partial \boldsymbol{K}^{\mathrm{Ele}}}{\partial \bar{\bvartheta}^{n+1}} \, \bar{\bPhi}_{m}^{n+1}
	+ \frac{\partial \boldsymbol{F}^{\mathrm{Ele}}}{\partial \bar{\bvartheta}^{n+1}} \right)_{\bar{\bPhi}_{m}^{n+1}, \bar{\bvartheta}_{m}^{n+1}, \bar{\bu}_{m}^{n+1}},
	\label{eq:Inexact_Linearization_Lagrangian_NL_Term_a_2} \\
	\displaystyle \left( \frac{\partial \bHH_1}{\partial \bar{\bu}^{n+1} } \right)_{\bar{\bv}_{m}^{n+1} } 
	&= \left( \frac{\partial \boldsymbol{K}^{\mathrm{Ele}}}{\partial \bar{\bu}^{n+1}} \, \bar{\bPhi}_{m}^{n+1} + \frac{\partial \boldsymbol{F}^{\mathrm{Ele}}}{\partial \bar{\bu}^{n+1}} \right)_{\bar{\bPhi}_{m}^{n+1}, \bar{\bvartheta}_{m}^{n+1}, \bar{\bu}_{m}^{n+1}},
	\label{eq:Inexact_Linearization_Lagrangian_NL_Term_a_3} \\
	\displaystyle \left( \frac{\partial \bHH_2}{\partial \bar{\bPhi}^{n+1} } \right)_{\bar{\bv}_{m}^{n+1} } 
	&= {\Delta t} \, \left( \frac{\partial \boldsymbol{F}^{\mathrm{The}}}{\partial \bar{\bPhi}^{n+1}} \right)_{\bar{\bPhi}_{m}^{n+1}, \bar{\bvartheta}_{m}^{n+1}, \bar{\bu}_{m}^{n+1}},
	\label{eq:Inexact_Linearization_Lagrangian_NL_Term_b_1} \\
	\displaystyle \left(\frac{\partial \bHH_2}{\partial \bar{\bvartheta}^{n+1} } \right)_{\bar{\bv}_{m}^{n+1} } 
	&= \left( \boldsymbol{M}^{\mathrm{The}} 
	+ {\Delta t} \left( \boldsymbol{K}^{\mathrm{The}} 
	+ \frac{\partial \boldsymbol{K}^{\mathrm{The}}}{\partial \bar{\bvartheta}^{n+1}} \bar{\bvartheta}_{m}^{n+1} 
	+ \frac{\partial \boldsymbol{F}^{\mathrm{The}}}{\partial \bar{\bvartheta}^{n+1}}
	\right) \right),
	\label{eq:Inexact_Linearization_Lagrangian_NL_Term_b_2} \\
	\displaystyle \left( \frac{\partial \bHH_2}{\partial \bar{\bu}^{n+1} } \right)_{\bar{\bv}_{m}^{n+1} } 
	&= {\Delta t} \left( \frac{\partial \boldsymbol{K}^{\mathrm{The}}}{\partial \bar{\bu}^{n+1}} \bar{\bvartheta}_{m}^{n+1} 
	+ \frac{\partial \boldsymbol{F}^{\mathrm{The}}}{\partial \bar{\bu}^{n+1}} \right)_{\bar{\bPhi}_{m}^{n+1}, \bar{\bvartheta}_{m}^{n+1}, \bar{\bu}_{m}^{n+1}},
	\label{eq:Inexact_Linearization_Lagrangian_NL_Term_b_3} \\
	\displaystyle \left( \frac{\partial \bHH_3}{\partial \bar{\bPhi}^{n+1} } \right)_{\bar{\bv}_{m}^{n+1} } 
	&= \left( \frac{\partial \boldsymbol{F}^{\mathrm{Mec}}}{\partial \bar{\bPhi}^{n+1}} \right)_{\bar{\bPhi}_{m}^{n+1}, \bar{\bu}_{m}^{n+1}},
	\label{eq:Inexact_Linearization_Lagrangian_NL_Term_c_1} \\
	\displaystyle \left(\frac{\partial \bHH_3}{\partial \bar{\bvartheta}^{n+1} } \right)_{\bar{\bv}_{m}^{n+1} } 
	&= \left( \frac{\partial \boldsymbol{K}^{\mathrm{Mec}}}{\partial \bar{\bvartheta}^{n+1}} \right)_{\bar{\bvartheta}_{m}^{n+1}, \bar{\bu}_{m}^{n+1}},
	\label{eq:Inexact_Linearization_Lagrangian_NL_Term_c_2} \\
	\displaystyle \left(\frac{\partial \bHH_3}{\partial \bar{\bu}^{n+1} } \right)_{\bar{\bv}_{m}^{n+1} } 
	&= \left( \frac{\partial \boldsymbol{K}^{\mathrm{Mec}}}{\partial \bar{\bu}^{n+1}} 
	+ \frac{\partial \boldsymbol{F}^{\mathrm{Mec}}}{\partial \bar{\bu}^{n+1}} \right)_{\bar{\bPhi}_{m}^{n+1}, \bar{\bvartheta}_{m}^{n+1}, \bar{\bu}_{m}^{n+1}}.
	\label{eq:Inexact_Linearization_Lagrangian_NL_Term_c_3}
\end{align}
In \eqref{eq:Inexact_Linearization_Lagrangian_NL_Term_c_3}, the first term is given by
\begin{equation}
	\left( \frac{\partial \boldsymbol{K}^{\mathrm{Mec}}}{\partial \bar{\bu}^{n+1}} \right)_{\bar{\bv}_{m}^{n+1} }= \sum_{e = 1}^{N^{\mathrm{Mec} } } {\bLL^{e}}^T 
	\left[ \int_{\Omega_e} (\bBB^{e})^T \, \left(\frac{\partial \bSSS_{EP}}{\partial \bEE} \colon \frac{\partial \bEE}{\partial \bar{\bu}^{n+1}} \right) \, \mathrm{d} \Omega_e \right].
\end{equation}
The Jacobian of the mechanical problem ${\partial \bSSS_{EP}}/{\partial \bEE}$ in \eqref{eq:Inexact_Linearization_Lagrangian_NL_Term_c_3} is also updated using the return mapping algorithm described in \cite{boatti-smp-16}. 
The pseudocode in Algorithm \ref{alg:smp_multiphysics} illustrates the flow of the numerical code used to solve the multiphysics problem.
\begin{center}
\begin{algorithm}[!]
\caption{Pseudocode of the multiphysics SMP problem}\label{alg:smp_multiphysics}
\begin{algorithmic}
\INPUT Mesh, current source $I_s(t)$ and constitutive laws.
\OUTPUT Mechanical, thermal and electromagnetic fields, and Joule losses
\Procedure{Multiphysics problem}{} 
    \State $t \gets t_{0}$, initialization of the physica fields $T|_{t_0} = T_0$ and $\bu|_{t_0} = \bu_0$,
    \For{$(n \gets 1$ To $N_{\mathrm{TS}} )$}              \Comment{\emph{the time loop $($index $n$$)$} \textcolor{white}{Innocent Ii}}
        \For{$(m \gets 1$ To $N_{\mathrm{NR}}^{\mathrm{M}} )$} \Comment{\emph{the NR loop of the overall problem}}
            \For{$(l \gets 1$ To $N_{\mathrm{RM}} )$}              \Comment{\emph{the NR loop of the return mapping}}
                \State Pass internal variables as input,
                \State Update the mechanical law $\bSS$ and ${\partial \bSSS_{EP}}/{\partial \bEE}$ using the RM
        \EndFor
        \State Update and compute the Jacobian, 
        \State Assemble $\left(\partial \bGG/\partial \bar{\bv}^{n+1} \right)_{\bar{\bv}_{m}^{n+1} }$ from \eqref{eq:Exact_Linearization_Lagrangian_NL_Term} and $\bGG(\bar{\bv}_{m}^{n+1})$ from \eqref{eq:Linearization_Lagrangian_NL_Term},
        \State Solve the system $\left(\partial \bGG/\partial \bar{\bv}^{n+1} \right)_{\bar{\bv}_{m}^{n+1} } \Delta \bv_{m+1}^{n+1} = -\bGG(\bar{\bv}_{m}^{n+1}, \bZZ_{L})$,
        \State compute the residual $\boldsymbol{r}(\bv_{m}^{n+1}) = \bGG(\bar{\bv}_{m}^{n+1})$,
        	\If{($||\boldsymbol{r}(\bv_{m}^{n+1})|| \leq \varepsilon_\mathrm{tol}$)}
    	    	\State Exit the nonlinear loop of the overall problem
      		\Else
        		\State Do another NR iteration for the overall problem
    		\EndIf
        \EndFor
    \State Go to the next time step $t \gets t + \Delta t$
    \EndFor
\EndProcedure
\end{algorithmic}
\end{algorithm}
\end{center}
%
%
%
%
%
%

\section{Numerical Tests}
\label{section:numerical_tests}

This section is devoted to the numerical testing of the electro-thermo-mechanical problem. A set of numerical tests similar to the ones developed in \cite{boatti-smp-16} are herein proposed. Whereas the authors in \cite{boatti-smp-16} considered ideal and non-ideal shape memory polymer materials, the main focus of this paper is on the fully coupled problem. Therefore, in Section \ref{sec:appli_mecha} we consider an ideal and a non-ideal shape memory polymer single element similar to the one in \cite{boatti-smp-16} for the validation of the thermomechanical problem. In Section \ref{sec:appli_fully_coupled} we only consider an ideal shape memory polymer stent for the fully coupled electro-thermo-mechanical problem.

The two considered tests are:
\begin{itemize}
    \item the uniaxial tests on a $1 \times 1 \times 1$ mm$^3$ single-element cube (SEC),
    \item the simulation of a cylindrical vascular stent (CVS) similar to the one described in \cite{boatti-smp-16}. 
    The stent has the same dimensions and material properties but without the small holes.
\end{itemize}
\begin{center}
    {
    \setlength\arrayrulewidth{1.0pt}
    \begin{table}[h]
        \begin{tabular}{|c|c|c|}
            \hline  
            Symbol          & Value & Unit\\
            \hline  
            \hline  
            $E^r$           & 0.9   & MPa \\  
            $E^g$           & 771   & MPa \\  
            $\nu^r$         & 0.49  & --  \\  
            $\nu^g$         & 0.29  & --  \\  
            $R^{pg}$        & 10    & MPa \\  
            $h$             & 0     & MPa \\  
            $\Delta \theta$ & 30 (SEC) -- 5 (CVS)     & K   \\  
            $\theta_t$      & 350 (SEC) -- 344 (CVS)  & K   \\  
            $w$             & 0.2 (SEC) -- 0.375 (CVS) & 1/K \\  
            \textcolor{white}{(ideal)} $c$ \textcolor{white}{(ideal)} & \textcolor{white}{(ideal)} 1 \textcolor{white}{(ideal)}   & \textcolor{white}{(ideal)} -- \textcolor{white}{(ideal)}  \\  
            $c^p$           & 0     & -- \\  
            \hline  
        \end{tabular}
        \caption{Model parameters of the mechanical problem}
        \label{table:mechanics}
    \end{table}
    }
\end{center}
Material properties listed in Table \ref{table:mechanics} are used for both cases. The software GetDP \cite{dular-getdp-98} was used to solve the fully coupled problem based on a total Lagrangian formulation.

\subsection{Validation of the mechanical problem}
\label{sec:appli_mecha}

Results of the thermo-mechanical model developed in \cite{boatti-smp-16} are reproduced. This model consisted of a temperature-dependent elasto-plastic model with a constant temperature field imposed for all Gauss points at any given time instant $t$. 
%
%

%
%
\pgfplotstableread{SE_I_Temperature_dat.txt}{\datatable}
\pgfplotstablecreatecol[copy column from table={SE_I_Temperature_dat.txt}{1}] {data_Time_Ideal} {\datatable}
\pgfplotstablecreatecol[copy column from table={SE_I_Temperature_dat.txt}{5}] {data_Temperature_Ideal} {\datatable}
\pgfplotstablecreatecol[copy column from table={SE_I_Displacement_x_dat.txt}{5}] {data_Displacement_x_Ideal} {\datatable}
\pgfplotstablecreatecol[copy column from table={SE_I_Displacement_y_dat.txt}{5}] {data_Displacement_y_Ideal} {\datatable}
\pgfplotstablecreatecol[copy column from table={SE_I_Displacement_z_dat.txt}{5}] {data_Displacement_z_Ideal} {\datatable}
\pgfplotstablecreatecol[copy column from table={SE_I_sigma_yy_dat.txt}{5}] {data_Sigma_Ideal} {\datatable}
\pgfplotstablecreatecol[copy column from table={SE_I_eps_small_yy_dat.txt}{5}] {data_GL_Ideal} {\datatable}
\pgfplotstablecreatecol[copy column from table={SE_I_PK1_yy_dat.txt}{5}] {data_PK1_Ideal} {\datatable}
\pgfplotstablecreatecol[copy column from table={SE_I_PK2_yy_dat.txt}{5}] {data_PK2_Ideal} {\datatable}
\pgfplotstablecreatecol[copy column from table={SE_NI_Temperature_dat.txt}{1}] {data_Time} {\datatable}
\pgfplotstablecreatecol[copy column from table={SE_NI_Temperature_dat.txt}{5}] {data_Temperature} {\datatable}
\pgfplotstablecreatecol[copy column from table={SE_NI_Displacement_x_dat.txt}{5}] {data_Displacement_x} {\datatable}
\pgfplotstablecreatecol[copy column from table={SE_NI_Displacement_y_dat.txt}{5}] {data_Displacement_y} {\datatable}
\pgfplotstablecreatecol[copy column from table={SE_NI_Displacement_z_dat.txt}{5}] {data_Displacement_z} {\datatable}
\pgfplotstablecreatecol[copy column from table={SE_NI_sigma_yy_dat.txt}{5}] {data_Sigma} {\datatable}
\pgfplotstablecreatecol[copy column from table={SE_NI_eps_small_yy_dat.txt}{5}] {data_GL} {\datatable}
\pgfplotstablecreatecol[copy column from table={SE_NI_PK1_yy_dat.txt}{5}] {data_PK1} {\datatable}
\pgfplotstablecreatecol[copy column from table={SE_NI_PK2_yy_dat.txt}{5}] {data_PK2} {\datatable}
\begin{figure}[]
	\begin{tikzpicture}
		\begin{axis}[xlabel={$\varepsilon_{yy}$ (\%)}, ylabel={$T$ (K)}, zlabel={$\sigma_{yy}$ (MPa)}, width=0.5\columnwidth]
			\addplot3[color=red, mark=cube*, mark size=2.0, mark options=solid] 		
					table[x expr=\thisrow{data_GL}*-1, y expr=\thisrow{data_Temperature}, z expr=\thisrow{data_Sigma}*1e-6]{\datatable};
			\addplot3[color=blue, mark=cube*, mark size=2.0, mark options=solid] 		
					table[x expr=\thisrow{data_GL_Ideal}*-1, y expr=\thisrow{data_Temperature_Ideal}, z expr=\thisrow{data_Sigma_Ideal}*1e-6]{\datatable};
		\end{axis}
	\end{tikzpicture}
	\begin{tikzpicture}
		\begin{axis}[xlabel={Time (s)}, ylabel={$\varepsilon_{yy}$ (\%)},width=0.5\columnwidth, legend style={at={(0.5,0.3)},anchor=north}]
			\addplot[color=red, mark=cube*, mark size=2.0, mark options=solid] 
					table[x expr=\thisrow{data_Time}, y expr=\thisrow{data_GL}*1]{\datatable}; 
			\addplot[color=blue, mark=cube*, mark size=2.0, mark options=solid] 
					table[x expr=\thisrow{data_Time_Ideal}, y expr=\thisrow{data_GL_Ideal}*1]{\datatable}; 
		\end{axis}
	\end{tikzpicture}
	\begin{tikzpicture}
		\begin{axis}[xlabel={Time (s)}, ylabel={$\sigma_{yy}$ (MPa)}, width=0.5\columnwidth]
			\addplot[color=red, mark=cube*, mark size=2.0, mark options=solid] 
					table[x expr=\thisrow{data_Time}, y expr=\thisrow{data_Sigma}*1e-6]{\datatable};
			\addplot[color=blue, mark=cube*, mark size=2.0, mark options=solid] 
					table[x expr=\thisrow{data_Time_Ideal}, y expr=\thisrow{data_Sigma_Ideal}*1e-6]{\datatable};
		\end{axis}
	\end{tikzpicture}
	\begin{tikzpicture}
		\begin{axis}[xlabel={$\varepsilon_{yy}$ (\%)}, ylabel={$\sigma_{yy}$ (MPa)}, width=0.5\columnwidth]
			\addplot[color=red, mark=cube*, mark size=2.0, mark options=solid] 
					table[x expr=\thisrow{data_GL}*1, y expr=\thisrow{data_Sigma}*1e-6]{\datatable};
			\addplot[color=blue, mark=cube*, mark size=2.0, mark options=solid] 
					table[x expr=\thisrow{data_GL_Ideal}*1, y expr=\thisrow{data_Sigma_Ideal}*1e-6]{\datatable};
		\end{axis}
	\end{tikzpicture}
	\begin{tikzpicture}
		\begin{axis}[xlabel={Temperature (s)}, ylabel={$\sigma_{yy}$ (MPa)}, width=0.5\columnwidth]
			\addplot[color=red, mark=cube*, mark size=2.0, mark options=solid] 
					table[x expr=\thisrow{data_Temperature}*1, y expr=\thisrow{data_Sigma}*1e-6]{\datatable};
			\addplot[color=blue, mark=cube*, mark size=2.0, mark options=solid] 
					table[x expr=\thisrow{data_Temperature_Ideal}*1, y expr=\thisrow{data_Sigma_Ideal}*1e-6]{\datatable};
		\end{axis}
	\end{tikzpicture}
	\begin{tikzpicture}
		\begin{axis}[xlabel={$T$ (K)}, ylabel={$\varepsilon_{yy}$ (\%)}, width=0.5\columnwidth]
			\addplot[color=red, mark=cube*, mark size=2.0, mark options=solid] 
					table[x expr=\thisrowno{5}*1, y expr=\thisrow{data_GL}*1]{\datatable}; 
			\addplot[color=blue, mark=cube*, mark size=2.0, mark options=solid] 
					table[x expr=\thisrowno{5}*1, y expr=\thisrow{data_GL_Ideal}*1]{\datatable}; 
		\end{axis}
	\end{tikzpicture}
    \caption{Results of the single element. 
    Top-left: temperature--strain--stress curve.
    Top-right: strain versus time curve.
    Middle-left: stress versus time curve.     
    Middle-right: strain versus stress curve. 
    Bottom-left: temperature versus stress curve.     
    Bottom-right: temperature versus strain curve.
    Blue curves correspond to the ideal case while red curves correspond to the non-ideal case.}
    \label{fig:smp-mechanics-single-element}
\end{figure}
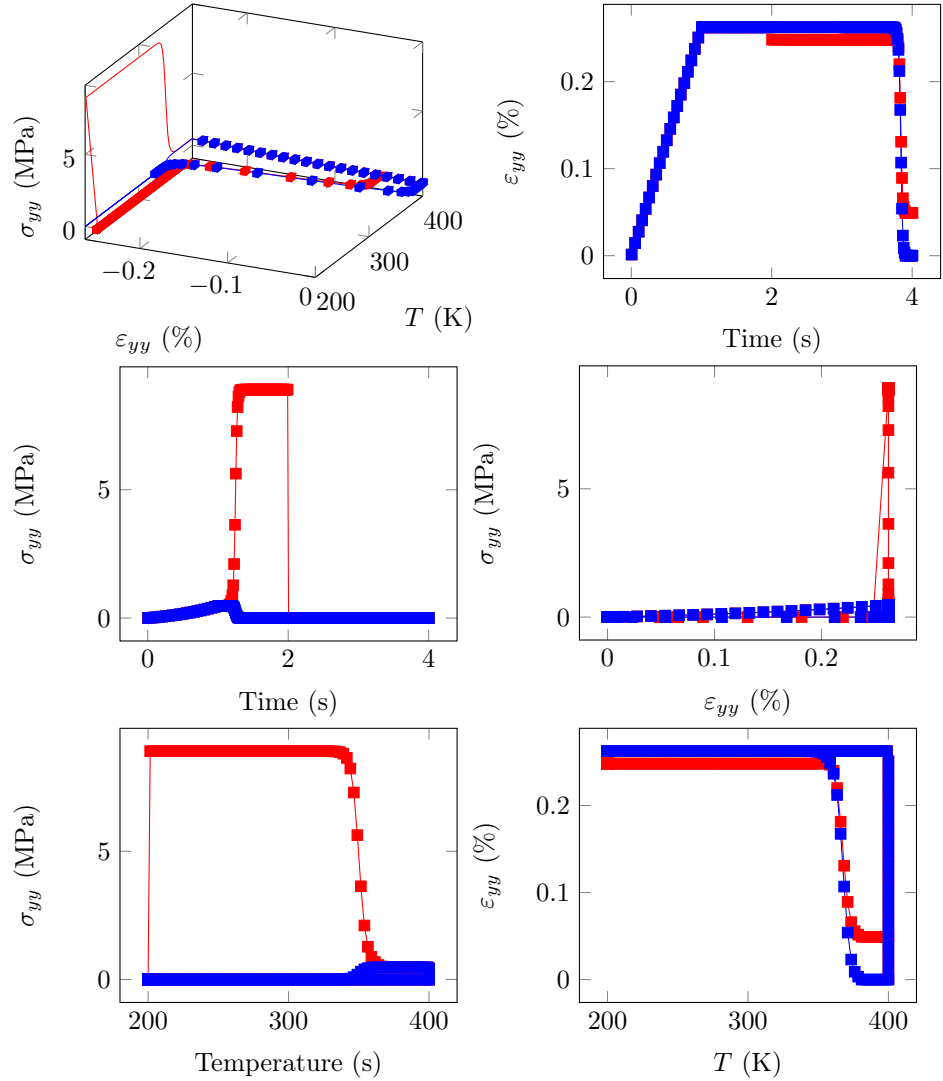 
%
%

%
%
Results of the single element are reported in Figure \ref{fig:smp-mechanics-single-element} for a high-temperature fixing similar to the  one used in \cite{boatti-smp-16}. The material is progressively deformed at 400K, then cooled down to 200K while keeping constant the deformation. The material is then unloaded at 200K before re-heating it up to 400K to trigger the shape-recovery (see Test 1 in Fig. 6 of \cite{boatti-smp-16}). The results reported in Figure \ref{fig:smp-mechanics-single-element} conform to those obtained in \cite{boatti-smp-16}.

\subsection{New results of the coupled problem}
\label{sec:appli_fully_coupled}

The results of the coupled problem are presented below. A description of the mechanical, thermal and electromagnetic problems is followed by the presentation of numerical results of the shape memory polymer. A best design can be obtained by choosing material properties for the thermal problem (thermal conductivity, mass density and heat capacity) that maintain a homogeneous temperature field in the stent, in order to avoid the appearance of regions with different phases during the recovery step. A non dimensionalization analysis of the thermal problem carried out in Section \ref{subsubsec:resultscoupledproblem} facilitates this design. However, the control of the temperature is complicated by the dependence of the mechanical stress on the temperature-dependent ratio of the glassy state $z^g(\theta_L)$.

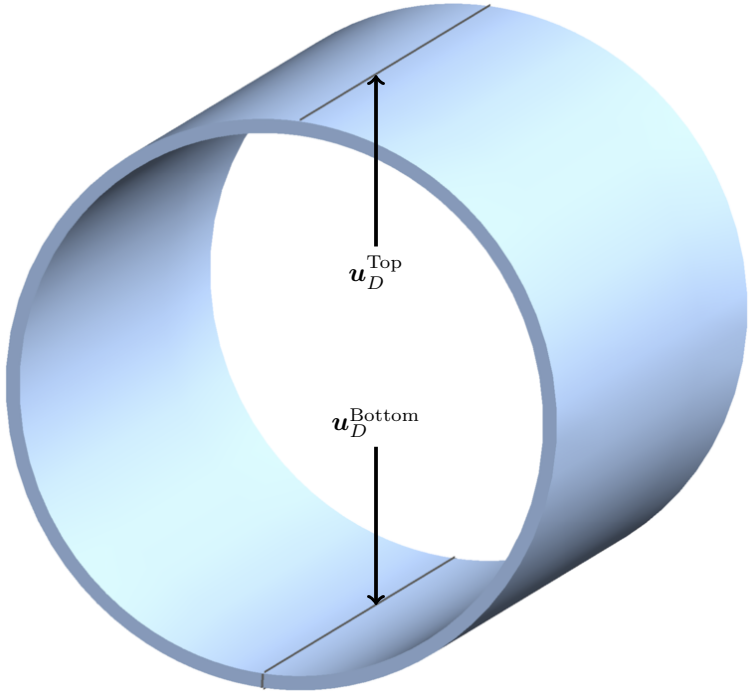
\begin{figure}[h]
    \centering
    \begin{tikzpicture}
        \draw (0, 0) node[inner sep=0] {\includegraphics[width=10cm]{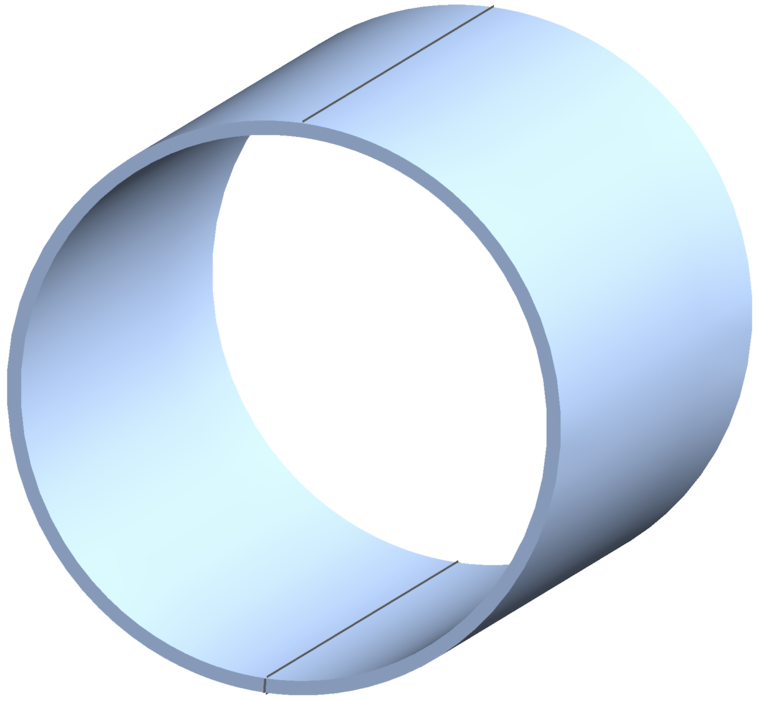}};
        \node (A) at (0.0, 1.0) {$\bu_D^{\mathrm{Top}}$};
        \node (B) at (0.0, 3.75) {};
        \node (C) at (-0.0, -1.0) {$\bu_D^{\mathrm{Bottom}}$};
        \node (D) at (-0.0, -3.55) {};
        \draw[->, line width=0.5mm, to path={-- (\tikztotarget)}](A) edge (B);
        \draw[->, line width=0.5mm, to path={-- (\tikztotarget)}](C) edge (D);
    \end{tikzpicture}
    \caption{Dirichlet boundary conditions for the mechanical problem. Zero displacement $\bu_D^{\mathrm{Bottom}}(t) = \boldmath{0}$ is imposed on the bottom section and a time-dependent displacement $\bu_D^{\mathrm{Top}}(t)$ similar to the one used for Test 1 of Fig. 6 in \cite{boatti-smp-16} is imposed on the top line.}
    \label{fig:smp-mechanics-stent}
\end{figure}
The geometry of Figure \ref{fig:smp-mechanics-stent} is used for the mechanical problem of the cylindrical vascular stent.

\subsubsection{Electromagnetic and thermal problems}

We simulate the insertion of a vascular shape memory polymer stent in a vein of the arm. The stent contains electric particles that can react to electromagnetic source fields produced by a coil wrapped around the arm by producing heat by the Joule effect. For the sake of simplicity, we consider the resulting shape memory polymer composite to be homogeneous with homogenized macroscopic material properties, thus ignoring the multiscale nature of the composite.

The mechanical problem is similar to the one in \cite{boatti-smp-16}, with the temperature field obtained by solving the thermal problem with the source generated by the eddy current losses. In the following, we define the electromagnetic and the thermal problems.

The temperature can be controlled by an electromagnetic field generated by a coil crossed by a current denoted $I_s(t)$. For all problems studied herein, we consider a single frequency source
\begin{equation}
	I_s(t) = I_0(t) (a + b \, \sin(\omega \, t)) = I_0(t) (a + b \, \sin(2 \, \pi \, f \, t)), 
\end{equation}  
where $I_0(t)$ (A) is piecewise, linear, time-dependent amplitude of the electric current, $\omega$ is the angular velocity and $f$ the frequency of the signal. The design parameters for the electromagnetic and thermal problems are the amplitude of the current, the frequency and the material properties: the electric conductivity $\sigma$ (S/m), the magnetic permeability $\mu$ (H/m), the mass density $\rho$, the heat capacity $c_p$ and the thermal conductivity $\bkappa$. In all our applications, we consider frequencies small than 1000Hz, $\sigma = 10^4$S/m and $\mu = \mu_0 \mu_{\mathrm{rel}}$ with $\mu_{\mathrm{rel}} = 20$, which corresponds to the wavelength $\lambda$ and skin depth $\delta$:
\begin{equation}
	\lambda = \frac{c}{f} = \frac{1}{\sqrt{\mu \epsilon} f } \approx 300\text{km} \quad , \quad \delta = \sqrt{\frac{2}{\mu \sigma \omega}} \approx 40\text{mm}. 
\end{equation}  
\begin{figure}[]
\centering
    \includegraphics[width=0.85\textwidth]{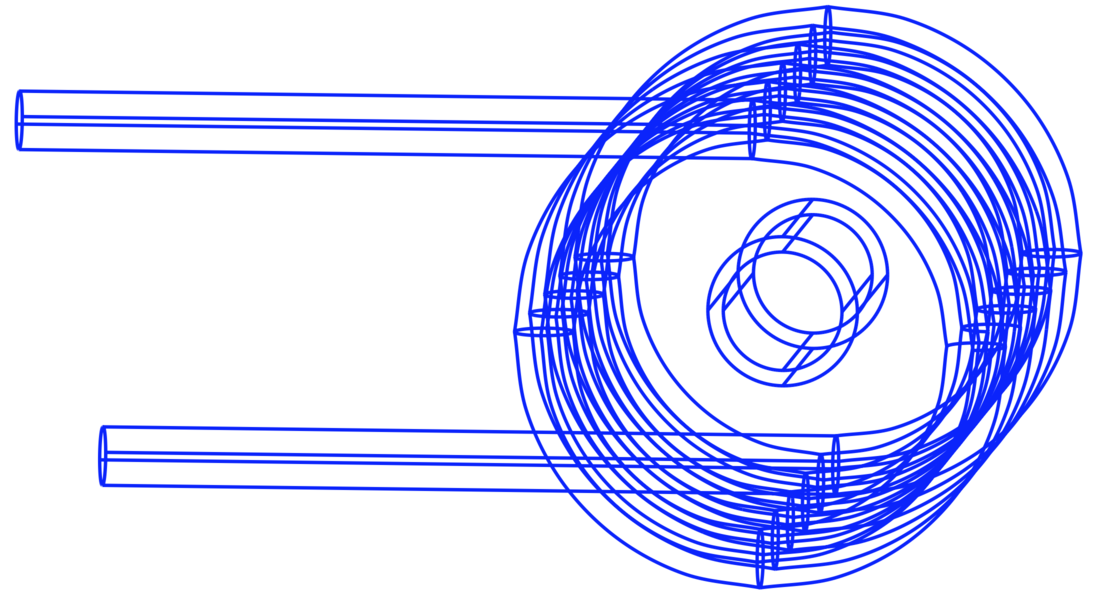}
    \includegraphics[width=0.85\textwidth]{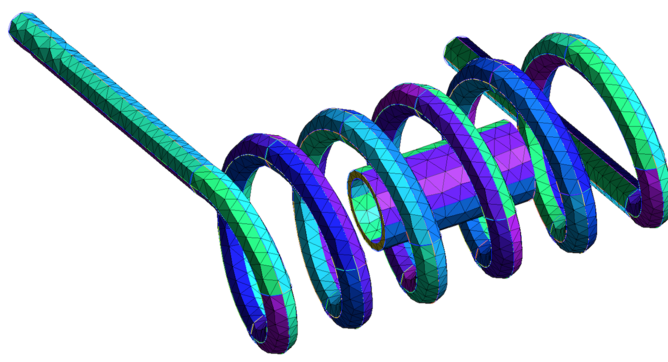}
    \includegraphics[width=0.75\textwidth]{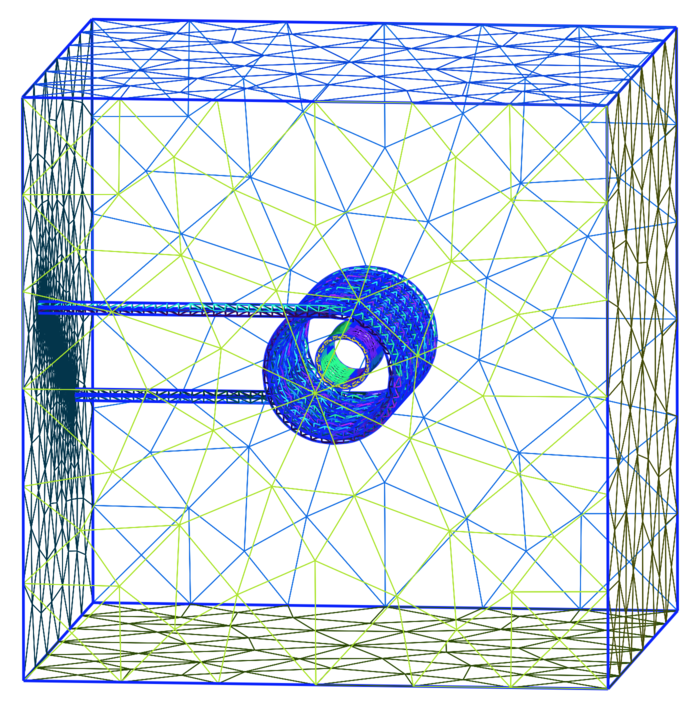}
    \caption{Geometry and mesh used for the coupled problem.
    Top: The cylindrical cardiovascular stent surrounded by an exciting coil.
    Middle: Mesh of the stent and the coil.
    Bottom: Mesh of the stent, the coil and the surrounding air. The enclosing box is used to bound the 
    computational domain for the electromagnetic problem assumed to be unbounded.}
    \label{fig:smp_physical_chemical}
\end{figure}
The wavelength is very large compared to the dimensions of the structure (typically 20mm for the length and 1mm for the thickness) that the \emph{quasistatic assumption} can be made \cite{rodriguez-10-eddycurrents}. Likewise, the skin depth is large compared to the dimensions of the stent that the eddy currents resulting from the reaction field can be neglected. Figure \ref{fig:smp_physical_chemical} illustrates the geometry used for the coupled problem.

To determine the magnetic induction source $\bb_s(t)$ for the electromagnetic problem, we consider a coil with a very large number of turns. The value of the magnetic field $\bh_s(t)$ and the magnetic induction $\bb_s(t)$ in the coil are homogeneous and given by \cite{field-coil}
\begin{equation}
	\bh_s(t) = h_s(t) \, \be_z = \frac{N}{L}I_s(t) \, \be_z \, \, , \, \, \bb_s(t) = b_s(t) \, \be_z = \mu \bh_s(t) = \mu \frac{N}{L}I_s(t) \, \be_z
    \label{eq:Applications_Sources_h_b}
\end{equation}
where $N$ is the number of turns, $L$ the length of the coil, $\mu$ the magnetic permeability of the material and $\be_z$ the direction oriented along the axis of the coil. From the Gauss magnetic law $\Div[]{\bb_s} = 0$, a source vector potential $\ba_s(t)$ can be derived from the magnetic induction $\bb_s(t)$ as $\bb_s(t) = \Curl[]{\ba_s(t)}$. In the computational domain of the stent, a possible vector potential that satisfies this equality and is symmetric with respect to the undeformed geometry of the stent is 
\begin{equation}
	\ba_s(x, y, t) = 0.5 \, b_s(t) (-y, x, 0) = 0.5 \, \mu \frac{N}{L} I_s(t) (-y, x, 0)
    \label{eq:Applications_Sources_a}
\end{equation}
with $x = X + u_x$ and $y = Y + u_y$ where $X$ and $Y$ are the coordinates expressed in the undeformed configuration and $u_x$ and $u_y$ are components of the displacement $\bu = (u_x, u_y, u_z)$. From \eqref{eq:Applications_Sources_h_b} and \eqref{eq:Applications_Sources_a}, it can be noted that  the magnetic field $\bh_s(t)$ and the magnetic induction $\bb_s(t)$ in the coil do not depend on spatial coordinates whereas the vector potential $\ba_s(t)$ depends on spatial coordinates $x$ and $y$. This vector potential is a one differential form that can be transformed as $\bAAA_s = \bFF^T \ba_s$ thus leading to the source in \eqref{eq:Simplified_Maxwell_Pullback_H_E}.

Table \ref{table:elemag_thermal} contains model parameters of the electromagnetic and thermal problems.
\begin{center}
    {
    \setlength\arrayrulewidth{1.0pt}
    \begin{table}[h]
        \begin{tabular}{|c|c|c|}
            \hline  
            Symbol      & Value    & Unit\\
            \hline  
            \hline  
            $I_0(t)$    & electric current waveform & A \\  
            $f$         & 1000                      & Hz \\  
            $\sigma$    & 10$^4$  & S/m \\  
            $\mu_r$     & 20    & --  \\  
            $\rho$      & 270   & kg m$^{-3}$ \\  
            $c_p$       & 10    & kg m$^2$ K$^{-1}$ s$^{-2}$ \\  
            \textcolor{white}{(ideal)} $k$ \textcolor{white}{(ideal)} & \textcolor{white}{(ideal)} 237 \textcolor{white}{(ideal)}   & \textcolor{white}{(ideal)} W m$^{-1}$ K$^{-1}$ \textcolor{white}{(ideal)}  \\  
            $h$        & 500    & W m$^{-2}$ K$^{-1}$ \\
            $N$        & 1000   & -- \\
            $L$        & 1      & m \\
            \hline  
        \end{tabular}
        \caption{Model parameters of the electromagnetic and thermal problems}
        \label{table:elemag_thermal}
    \end{table}
    }
\end{center}

Defining the thermal problem resulting from the deployment of the actual stent is challenging. Though it is easy to control the temperature of the device during the first three stages (\emph{loading} at high temperature, \emph{cooling} followed by \emph{unloading/insertion} of the stent) most of which are done outside the human body, the last step, the \emph{recovery}, necessitates controlling the temperature using electromagnetic fields. In this paper, we simulate the control of the entire deployment process using the electromagnetic fields.

During the last step, different modes of heat exchange can be considered: (1) heat conduction in the stent, at the interface of the stent and the surrounding tissue and in the tissue itself and (2) forced convection at part of the boundary of surface of the stent in contact with the blood flowing in the vein. The surface of the stent in contact with the tissue/blood varies during the process of recovery and its detection would necessitate consideration of contact mechanics. For the sake of simplicity, we only consider forced convection.

Finally, thanks to the assumptions made of the electromagnetic and thermal problems, all three problems can only be solved on the computational domain of the stent thus neglecting the surrounding environment.

\subsubsection{Results of the coupled problems}
\label{subsubsec:resultscoupledproblem}

Results of the coupled problem are herein reported. As mentioned earlier, the main difference between this section and section \ref{sec:appli_mecha} lies in the use of a temperature field obtained by solving the heat equation on a moving domain with the source obtained by solving the electromagnetic problem instead of a priori imposing a temperature field at each time instant $t$.

\begin{figure}[!]
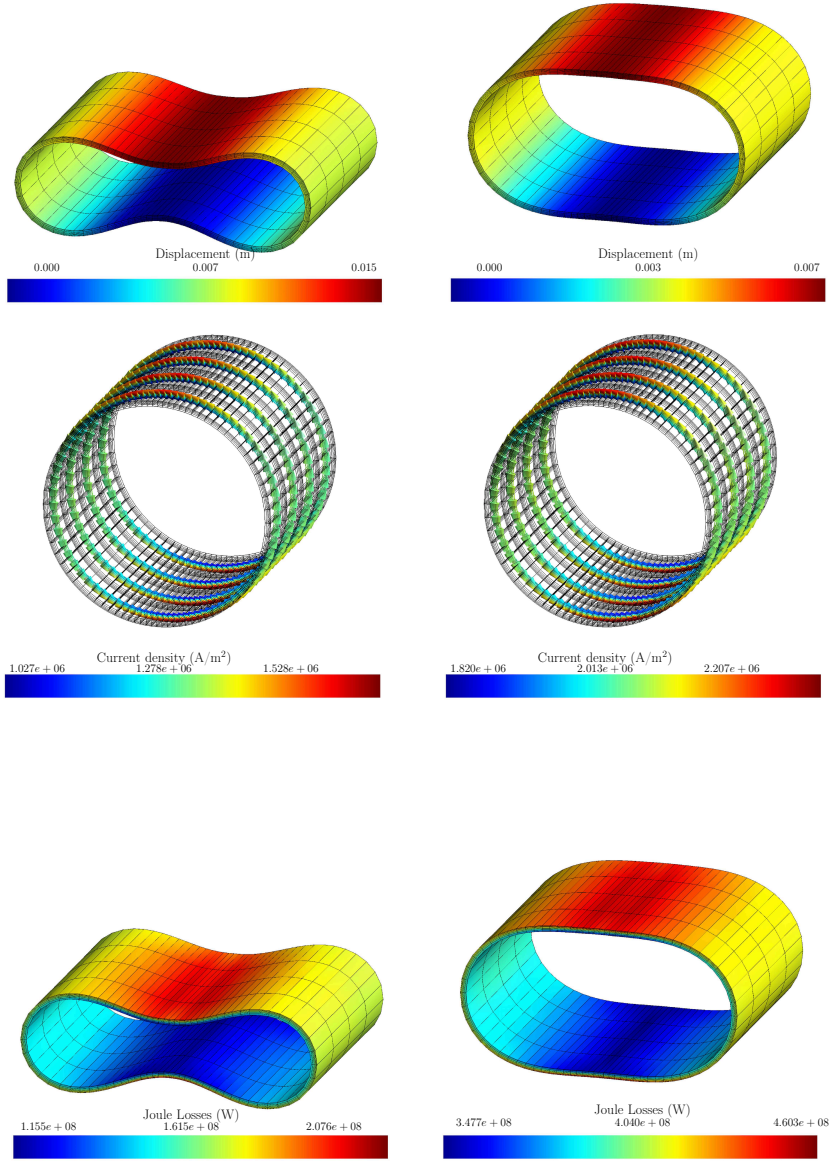

    \begin{center}
    	\scalebox{0.30}{\input{u_TS_89.tex}}
    	\scalebox{0.30}{\input{u_TS_91.tex}}
    \end{center}
    \vspace{10mm}
    \begin{center}
    	\hspace{-12mm}\scalebox{0.30}{\input{J_Tot_TS_89.tex}}
    	\hspace{10mm}
    	\scalebox{0.30}{\input{J_Tot_TS_91.tex}}
    \end{center}
    \vspace{10mm}
    \begin{center}
    	\scalebox{0.30}{\input{Joule_Losses_TS_89.tex}}
    	\hspace{3mm}
    	\scalebox{0.30}{\input{Joule_Losses_TS_91.tex}}
    \end{center}
	\caption{Physical fields during the recovery. The left images correspond to $t = 4.78125 \times 10^{-3}$s and the right images correspond to $t = 4.84375 \times 10^{-3}$s. 
		Top : displacement $\bu$, middle : current density $\bJJJ$, bottom : Joule losses $w_L$.
		}
\label{fig:Applications_Fields_bu_bj}
\end{figure}

\begin{figure}[!]
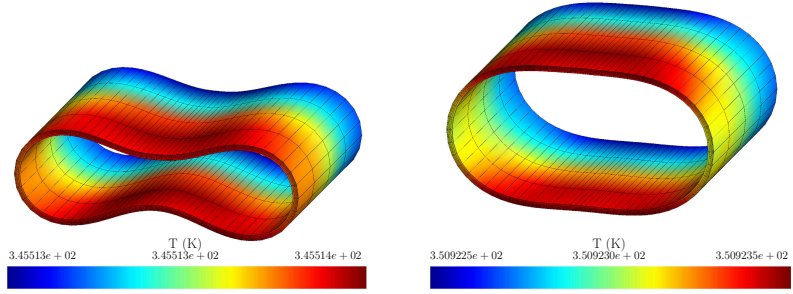

    \begin{center}
    	\scalebox{0.30}{\input{T_TS_89.tex}}
    	\hspace{2mm}
    	\scalebox{0.30}{\input{T_TS_91.tex}}
    \end{center}
	\caption{Temperature at $t = 4.78125 \times 10^{-3}$s (left) and $t = 4.84375 \times 10^{-3}$s (right).
		}
\label{fig:Applications_Fields_JL_T}
\end{figure}

Figures \ref{fig:Applications_Fields_bu_bj} and \ref{fig:Applications_Fields_JL_T} show the displacement $\bu$, the current density $\bj$, Joule losses and the temperature $T$ at the instances $t = 4.78125 \times 10^{-3}$s and $t = 4.84375 \times 10^{-3}$s. 
This can play an important role in the design of the stent, especially for the computation of the temperature. Indeed, the high dependency of the ratio of the glassy states on the temperature $z^g(\theta_L)$ necessitates selecting electromagnetic and thermal loadings as well as thermal material properties that allow for a quick diffusion of the heat sources throughout the stent, to avoid inhomogeneities of temperature that would cause different regions of the stent to be in different phases (rubbery/glassy) during the recovery step. 
Another issue concerns the use of the Newton--Raphson method to solve the nonlinear coupled problem. Considerably large and inhomogeneous increments of temperature computed especially during the first nonlinear iterations of the Newton--Raphson scheme may lead to inhomogeneities of temperature and cause slow convergence in the recovery process.

To avoid inhomogeneities of temperature in the recovery step, we developed the following normalization process, which makes the problem well conditioned.

The process starts with the linearized version of the \emph{heat equation} \eqref{eq:Strong_Heat_Governing_Equation_Lagrangian}:
\begin{equation}
    \displaystyle \rho_L c_p \frac{\partial \theta_L}{\partial t} + \Divm[_X]{ \left[ \bkappa_L 
    \Gradm[_X]{\theta_L} \right] } = -w_L(\phi, \theta_L, \bu).
	\label{eq:Weak_Thermal_Equation_Lagrangian_Dimensionless_Eq_1}
\end{equation}
A new coordinate system ($\tau$, $\boldsymbol{\eta}$) is introduced as:
\begin{multline}
    t = T_{c} \, \tau \, \, , \, \, dt = T_{c} \, d \tau \, \, , \, \,  
    \frac{\partial (\cdot)}{\partial t} = \frac{1}{T_{c}} \frac{\partial (\cdot)}{\partial \tau}, \\ 
    X_i = L_{c} \, \eta_i  \, \, , \, \,  dX_i = L_{c} \, d \eta_i \, \, , \, \, 
    \frac{\partial (\cdot)}{\partial X_i} = \frac{1}{L_{c}} \frac{\partial (\cdot)}{\partial \eta_i}, \quad \quad 
    \theta_L = \theta_{c} \, \bar{\theta}_L
    \label{eq:Weak_Thermal_Equation_Lagrangian_Dimensionless_Eq_2}
\end{multline}
where $\bar{\theta}$, $T_{c}$, $L_{c}$, $\tau$, and $\eta_i$ are the characteristic temperature, the characteristic time, the characteristic length, the dimensionless temporal and spatial coordinates, respectively. The derivatives in \eqref{eq:Weak_Thermal_Equation_Lagrangian_Dimensionless_Eq_1} are transformed as:
\begin{multline}
    \displaystyle \frac{\partial \theta_L}{\partial t} = \frac{1}{T_{c}} \frac{\partial \left(\theta_{c} \, 
    \bar{\theta}_L \right) }{\partial \tau} = \frac{\theta_{c}}{T_{c}} \frac{\partial \bar{\theta}_L}{\partial \tau} \quad , \quad 
    \frac{\partial \theta_L}{\partial X_i} = \frac{1}{L_{c}} 
    \frac{\partial \left(\theta_{c} \, \bar{\theta}_L \right) }{\partial \eta_i} = \frac{\theta_{c}}{L_{c}} 
    \frac{\partial \bar{\theta}_L}{\partial \eta_i}, \\
    \displaystyle \frac{\partial^2 \theta_L}{\partial X_i^2} = \frac{1}{L_{c}} 
    \frac{\partial \left( \displaystyle \frac{\partial \theta_L}{\partial x_i} \right) }{\partial \eta_i} = 
    \frac{\theta_{c}}{L_{c}^2} \frac{\partial^2 \bar{\theta}_L}{\partial \eta_i^2}
\label{eq:Weak_Thermal_Equation_Lagrangian_Dimensionless_Eq_3}
\end{multline}
which leads to the \emph{dimensionless heat equation}:
\begin{equation}
    \displaystyle \frac{\partial \bar{\theta}_L}{\partial \tau} + 
    \frac{T_{c} \bar{\bkappa}_L}{\rho_L c_p L_{c}^2} \Divm[_{\eta}]{ \left[ 
    \Gradm[_{\eta}]{\bar{\theta}_L} \right] } = - \frac{T_{c}}{\rho_L c_p L_{c}^2} \bar{w}_L(\bar{\phi}, \bar{\theta}_L, \bar{\bu})
\label{eq:Weak_Thermal_Equation_Lagrangian_Dimensionless_Eq_4}
\end{equation}
where the thermal conductivity was assumed constant and barred quantities are defined in the new coordinate system. For the first two terms to be of the same order of magnitude, i.e., for the temperature to have enough time to diffuse in the stent, the material properties must be chosen such that
\begin{equation}
    T_{c} = \displaystyle \frac{\rho_L c_p L_{c}^2}{\bar{\bkappa}_L} .
\label{eq:Weak_Thermal_Equation_Lagrangian_Dimensionless_Eq_5}
\end{equation}
%
%

%
%
%
%
%
\pgfplotstableread{Stent_Coupled_Temperature_dat.txt}{\datatable}
\pgfplotstablecreatecol[copy column from table={Stent_Coupled_Temperature_dat.txt}{5}] {data_Temperature} {\datatable}
\pgfplotstablecreatecol[copy column from table={Stent_Coupled_Temperature_dat.txt}{1}] {data_Time} {\datatable}
\pgfplotstablecreatecol[copy column from table={Stent_Coupled_Displacement_x_dat.txt}{5}] {data_Displacement_x} {\datatable}
\pgfplotstablecreatecol[copy column from table={Stent_Coupled_Displacement_y_dat.txt}{5}] {data_Displacement_y} {\datatable}
\pgfplotstablecreatecol[copy column from table={Stent_Coupled_Displacement_z_dat.txt}{5}] {data_Displacement_z} {\datatable}
\pgfplotstablecreatecol[copy column from table={Stent_Coupled_Reaction_Force_1_dat.dat}{1}] {reaction_force_1} {\datatable}
\pgfplotstablecreatecol[copy column from table={Stent_Coupled_Reaction_Force_2_dat.dat}{1}] {reaction_force_2} {\datatable}
\pgfplotstablecreatecol[copy column from table={Stent_Coupled_Reaction_Force_3_dat.dat}{1}] {reaction_force_3} {\datatable}
\pgfplotstablecreatecol[copy column from table={Stent_Coupled_sigma_yy_dat.txt}{5}] {data_Sigma} {\datatable}
\pgfplotstablecreatecol[copy column from table={Stent_Coupled_GL_yy_dat.txt}{5}] {data_GL} {\datatable}
\pgfplotstablecreatecol[copy column from table={Stent_Coupled_PK2_yy_dat.txt}{5}] {data_PK2} {\datatable}
\pgfplotstableread{Stent_NonCoupled_Temperature_dat.txt}{\datatableMech}
\pgfplotstablecreatecol[copy column from table={Stent_NonCoupled_Temperature_dat.txt}{5}] {data_Temperature_Mech} {\datatableMech}
\pgfplotstablecreatecol[copy column from table={Stent_NonCoupled_Temperature_dat.txt}{1}] {data_Time_Mech} {\datatableMech}
\pgfplotstablecreatecol[copy column from table={Stent_NonCoupled_Displacement_x_dat.txt}{5}] {data_Displacement_x_Mech} {\datatableMech}
\pgfplotstablecreatecol[copy column from table={Stent_NonCoupled_Displacement_y_dat.txt}{5}] {data_Displacement_y_Mech} {\datatableMech}
\pgfplotstablecreatecol[copy column from table={Stent_NonCoupled_Displacement_z_dat.txt}{5}] {data_Displacement_z_Mech} {\datatableMech}
\pgfplotstablecreatecol[copy column from table={Stent_NonCoupled_Reaction_Force_1_dat.dat}{1}] {reaction_force_1_Mech} {\datatableMech}
\pgfplotstablecreatecol[copy column from table={Stent_NonCoupled_Reaction_Force_2_dat.dat}{1}] {reaction_force_2_Mech} {\datatableMech}
\pgfplotstablecreatecol[copy column from table={Stent_NonCoupled_Reaction_Force_3_dat.dat}{1}] {reaction_force_3_Mech} {\datatableMech}
\pgfplotstablecreatecol[copy column from table={Stent_NonCoupled_sigma_yy_dat.txt}{5}] {data_Sigma_Mech} {\datatable}
\pgfplotstablecreatecol[copy column from table={Stent_NonCoupled_GL_yy_dat.txt}{5}] {data_GL_Mech} {\datatable}
\pgfplotstablecreatecol[copy column from table={Stent_NonCoupled_PK1_yy_dat.txt}{5}] {data_PK1_Mech} {\datatable}
\pgfplotstablecreatecol[copy column from table={Stent_NonCoupled_PK2_yy_dat.txt}{5}] {data_PK2_Mech} {\datatable}
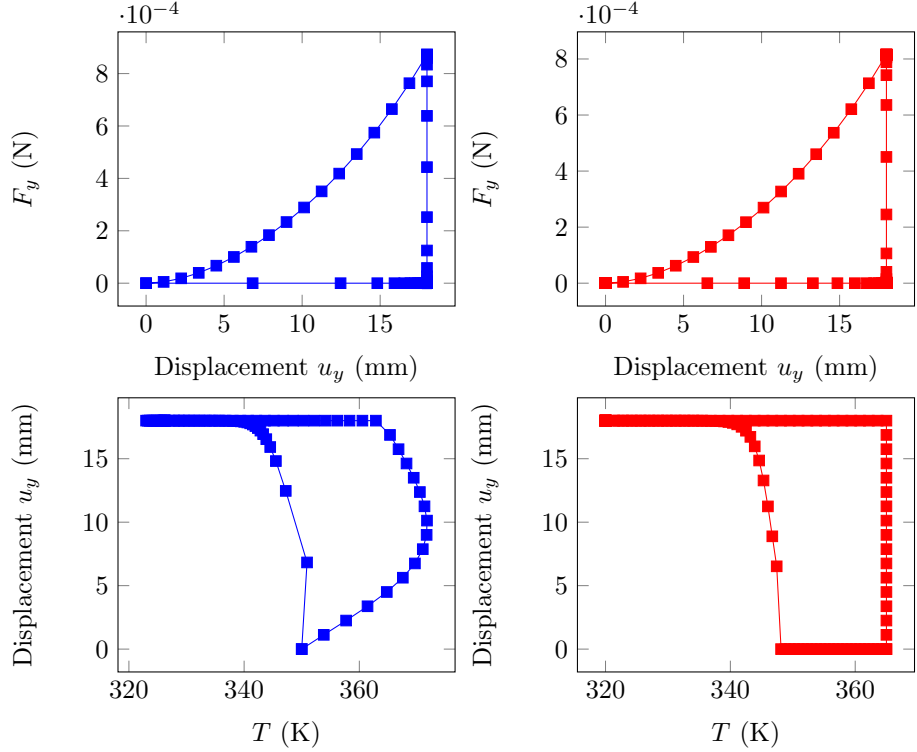
\begin{figure}[!]
	\begin{tikzpicture}
		\begin{axis}[xlabel={Displacement $u_y$ (mm)}, ylabel={$F_y$ (N)},width=0.5\columnwidth]
			\addplot[color=blue, mark=cube*, mark size=2.0, mark options=solid] 
					table[x expr=\thisrow{data_Displacement_y}*-1e3, y expr=\thisrow{reaction_force_2}*1]{\datatable};
		\end{axis}
	\end{tikzpicture}
	\begin{tikzpicture}
		\begin{axis}[xlabel={Displacement $u_y$ (mm)}, ylabel={$F_y$ (N)},width=0.5\columnwidth]
			\addplot[color=red, mark=cube*, mark size=2.0, mark options=solid] 
					table[x expr=\thisrow{data_Displacement_y_Mech}*-1e3, y expr=\thisrow{reaction_force_2_Mech}*1]{\datatableMech};
		\end{axis}
	\end{tikzpicture}
	\begin{tikzpicture}
		\begin{axis}[xlabel={$T$ (K)}, ylabel={Displacement $u_y$ (mm)},width=0.5\columnwidth]
			\addplot[color=blue, mark=cube*, mark size=2.0, mark options=solid] 
					table[x expr=\thisrowno{5}*1, y expr=\thisrow{data_Displacement_y}*-1e3]{\datatable};
		\end{axis}
	\end{tikzpicture}
	\begin{tikzpicture}
		\begin{axis}[xlabel={$T$ (K)}, ylabel={Displacement $u_y$ (mm)},width=0.5\columnwidth]
			\addplot[color=red, mark=cube*, mark size=2.0, mark options=solid] 
					table[x expr=\thisrow{data_Temperature_Mech}, y expr=\thisrow{data_Displacement_y_Mech}*-1e3]{\datatableMech};
		\end{axis}
	\end{tikzpicture}
	\caption{Results of the coupled problem. 
	Top : Force versus displacement curves. Computed temperature (left) and imposed temperature (right). 
	Bottom: Displacement versus temperature curves. Computed temperature (left) and imposed temperature (right).}
    \label{fig:smp-coupled-stent_1}
\end{figure} 
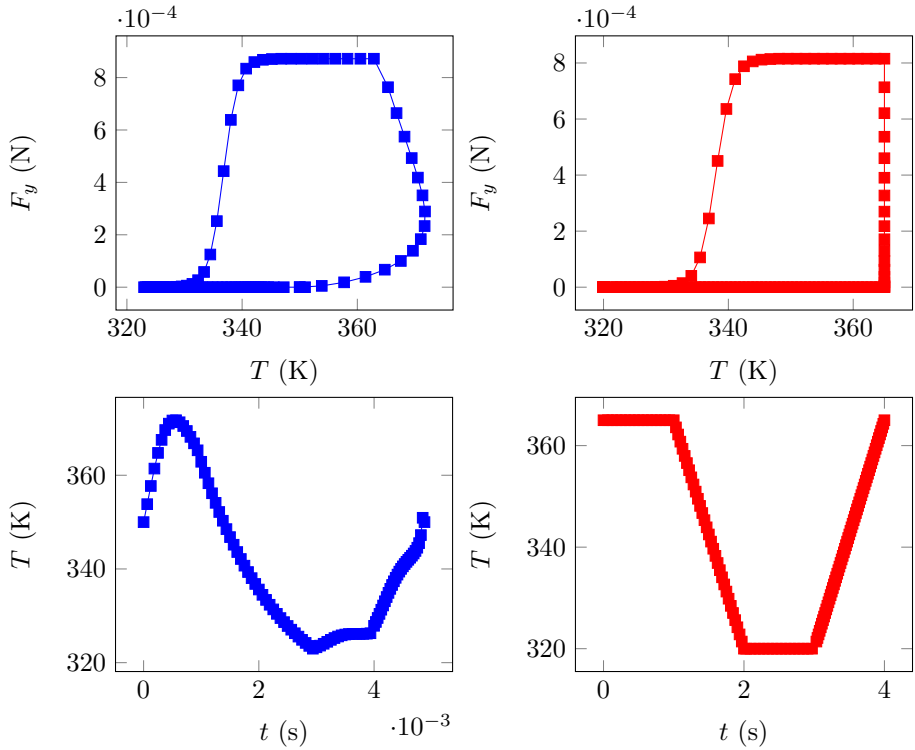
\begin{figure}[h]
	\begin{tikzpicture}
		\begin{axis}[xlabel={$T$ (K)}, ylabel={$F_y$ (N)},width=0.5\columnwidth]
			\addplot[color=blue, mark=cube*, mark size=2.0, mark options=solid] 
					table[x expr=\thisrow{data_Temperature}, y expr=\thisrow{reaction_force_2}*1]{\datatable};
		\end{axis}
	\end{tikzpicture}
	\begin{tikzpicture}
		\begin{axis}[xlabel={$T$ (K)}, ylabel={$F_y$ (N)},width=0.5\columnwidth]
			\addplot[color=red, mark=cube*, mark size=2.0, mark options=solid] 
					table[x expr=\thisrow{data_Temperature_Mech}, y expr=\thisrow{reaction_force_2_Mech}*1]{\datatableMech};
		\end{axis}
	\end{tikzpicture}
	\begin{tikzpicture}
		\begin{axis}[xlabel={$t$ (s)}, ylabel={$T$ (K)},width=0.5\columnwidth]
			\addplot[color=blue, mark=cube*, mark size=2.0, mark options=solid] 
					table[x expr=\thisrow{data_Time}, y expr=\thisrow{data_Temperature}]{\datatable};
		\end{axis}
	\end{tikzpicture}
	\begin{tikzpicture}
		\begin{axis}[xlabel={$t$ (s)}, ylabel={$T$ (K)},width=0.5\columnwidth]
			\addplot[color=red, mark=cube*, mark size=2.0, mark options=solid] 
					table[x expr=\thisrow{data_Time_Mech}, y expr=\thisrow{data_Temperature_Mech}]{\datatableMech};
		\end{axis}
	\end{tikzpicture}
	\caption{Results of the coupled problem. 
	Top : Force versus temperature curves. Computed temperature (left) and imposed temperature (right).
	Bottom : Temperature versus time curves. Computed temperature (left) and imposed temperature (right).}
    \label{fig:smp-coupled-stent_2}
\end{figure} 
%
%
%
%
Results of the coupled problem are reported in Figures \ref{fig:smp-coupled-stent_1}-\ref{fig:smp-coupled-stent_2} for a stent with the high-temperature fixing and slightly different material properties as those reported in \cite{boatti-smp-16}. In the case of the imposed temperature, the material is progressively deformed at 350K, then cooled down to 320K while the deformation is maintained constant, then unloaded at 320K, and finally re-heated up to 350K to trigger shape-recovery. We mimick the same trajectory of the temperature by changing material properties, the frequency and the amplitude of the excitation source. Results of the coupled problem are different from the ones obtained with the imposed temperature. An optimal control of the temperature using the source current $I_s(t)$ and the geometry of the coil as control parameters can allow to prescribe a temperature trajectory convenient for surgical purposes.
\clearpage 
%
%
\section{Conclusions}
\label{section:conclusions}
%
%
In this paper, the deployment of a vascular shape memory polymer stent in a vein of an arm is simulated. The temperature field used in the thermo-mechanical model of the stent is controlled by solving for electromagnetic fields generated by a coil wrapped around the arm.  
The controllability of the temperature depends on the choice of the electromagnetic source field determined by the amplitude and the excitation frequency of the current flowing through the coil, and on the material properties used for the thermal problem. An initial design of the stent which allows for the diffusion of heat and leads to a homogeneous distribution of temperature during the recovery step is proposed. 
The optimal control of the temperature of the devices can further be carried out, thus allowing the device to follow a prescribed temperature trajectory that might be convenient for surgical purposes.
%
%
%
%
\section*{Acknowledgments}
\label{sec:motivation}
The authors would like to thank Dr. Elisa Boatti at Georgia Institute of Technology, Prof. Ludovic Noels and Miguel Pareja Mu$\tilde{\text{n}}$oz at the University of Li\`{e}ge for fruitful discussions on mechanical models of shape memory polymers. They would also like to thank Prof. Christophe Geuzaine at the University of Li\`{e}ge for discussions about the implementation in GetDP.
The first author is particularly indebted to Dr. Fran\c{c}ois Henrotte at the University of Li\`{e}ge for the discussions on electromagnetic formulations under large deformations. 
During the time the research was carried out, Innocent Niyonzima was a postdoctoral Fellow with the Belgian American Educational Foundation (BAEF). He is also partially supported by an excellence grant from Wallonie-Bruxelles International (WBI).
%
%
\bibliography{cmame_smp}

\newcommand{\noop}[1]{}
\begin{thebibliography}{10}
\expandafter\ifx\csname url\endcsname\relax
  \def\url#1{\texttt{#1}}\fi
\expandafter\ifx\csname urlprefix\endcsname\relax\def\urlprefix{URL }\fi
\expandafter\ifx\csname href\endcsname\relax
  \def\href#1#2{#2} \def\path#1{#1}\fi

\bibitem{uweb-biometarials-04}
An introduction to biomaterials,
  \url{https://www.uweb.engr.washington.edu/research/tutorials/introbiomat.html},
  accessed: 2017-12-22.

\bibitem{bhatia-biomaterials-10}
S.~K. Bhatia, Biomaterials for clinical applications, Springer Science \&
  Business Media, 2010.

\bibitem{rezaie-biomaterials-15}
H.~R. Rezaie, L.~Bakhtiari, A.~{\"O}chsner, Biomaterials and their
  applications, Springer, 2015.

\bibitem{huang-smp-11}
W.~M. Huang, B.~Yang, Y.~Q. Fu, Polyurethane shape memory polymers, CRC Press,
  2011.

\bibitem{small-biomaterials-07}
W.~Small, P.~R. Buckley, T.~S. Wilson, W.~J. Benett, J.~Hartman, D.~Saloner,
  D.~J. Maitland, Shape memory polymer stent with expandable foam: a new
  concept for endovascular embolization of fusiform aneurysms, IEEE
  Transactions on Biomedical Engineering 54~(6) (2007) 1157--1160.

\bibitem{maitland-biomaterials-07}
D.~J. Maitland, W.~Small, J.~M. Ortega, P.~R. Buckley, J.~Rodriguez,
  J.~Hartman, T.~S. Wilson, Prototype laser-activated shape memory polymer foam
  device for embolic treatment of aneurysms, Journal of biomedical optics
  12~(3) (2007) 030504--030504.

\bibitem{yakacki-stents-07}
C.~M. Yakacki, R.~Shandas, C.~Lanning, B.~Rech, A.~Eckstein, K.~Gall,
  Unconstrained recovery characterization of shape-memory polymer networks for
  cardiovascular applications, Biomaterials 28~(14) (2007) 2255--2263.

\bibitem{baer-stents-09}
G.~M. Baer, T.~S. Wilson, W.~Small, J.~Hartman, W.~J. Benett, D.~L. Matthews,
  D.~J. Maitland, Thermomechanical properties, collapse pressure, and expansion
  of shape memory polymer neurovascular stent prototypes, Journal of Biomedical
  Materials Research Part B: Applied Biomaterials 90~(1) (2009) 421--429.

\bibitem{ajili-stents_09}
S.~H. Ajili, N.~G. Ebrahimi, M.~Soleimani, Polyurethane/polycaprolactane blend
  with shape memory effect as a proposed material for cardiovascular implants,
  Acta biomaterialia 5~(5) (2009) 1519--1530.

\bibitem{yahia-smp-15}
L.~Yahia, Shape memory polymers for biomedical applications, Elsevier, 2015.

\bibitem{wischke-drugdelivery-10}
C.~Wischke, A.~T. Neffe, S.~Steuer, A.~Lendlein, Comparing techniques for drug
  loading of shape-memory polymer networks--effect on their functionalities,
  European Journal of Pharmaceutical Sciences 41~(1) (2010) 136--147.

\bibitem{nagahama-biomaterial-09}
K.~Nagahama, Y.~Ueda, T.~Ouchi, Y.~Ohya, Biodegradable shape-memory polymers
  exhibiting sharp thermal transitions and controlled drug release,
  Biomacromolecules 10~(7) (2009) 1789--1794.

\bibitem{fish-piezo-03}
J.~Fish, W.~Chen, Modeling and simulation of piezocomposites, Computer methods
  in applied mechanics and engineering 192~(28-30) (2003) 3211--3232.

\bibitem{elhadrouz-piezo-06}
M.~Elhadrouz, T.~B. Zineb, E.~Patoor, Finite element analysis of a multilayer
  piezoelectric actuator taking into account the ferroelectric and ferroelastic
  behaviors, International journal of engineering science 44~(15) (2006)
  996--1006.

\bibitem{anderson-coupling-07}
P.~I. Anderson, A.~J. Moses, H.~J. Stanbury, Assessment of the stress
  sensitivity of magnetostriction in grain-oriented silicon steel, IEEE
  transactions on magnetics 43~(8) (2007) 3467--3476.

\bibitem{khalaquzzaman-piezoelectric-multiscale-12}
M.~Khalaquzzaman, B.~Xu, S.~Ricker, R.~M{\"u}ller, {Computational
  homogenization of piezoelectric materials using FE$^2$ to determine
  configurational forces}, Technische Mechanik 32~(1) (2012) 21--37.

\bibitem{kuznetsov-hmm-12}
S.~Kuznetsov, J.~Fish, Mathematical homogenization theory for electroactive
  continuum, International Journal for Numerical Methods in Engineering 91~(11)
  (2012) 1199--1226.

\bibitem{perevertov-coupling-15}
O.~Perevertov, J.~Thielsch, R.~Sch{\"a}fer, {Effect of applied tensile stress
  on the hysteresis curve and magnetic domain structure of grain-oriented
  transverse Fe-3\% Si steel}, Journal of Magnetism and Magnetic Materials 385
  (2015) 358--367.

\bibitem{bishay-piezo-electro-magnetic-15}
P.~L. Bishay, S.~N. Atluri, {Computational Piezo-Grains (CPGs) for a
  highly-efficient micromechanical modeling of heterogeneous
  piezoelectric--piezomagnetic composites}, European Journal of
  Mechanics-A/Solids 53 (2015) 311--328.

\bibitem{eringen-electrodynamics-12}
A.~C. Eringen, G.~A. Maugin, Electrodynamics of continua I: foundations and
  solid media, Springer Science \& Business Media, 2012.

\bibitem{pao-electrodynamics-78}
Y.-H. Pao, Electromagnetic forces in deformable continua, in: Mechanics today,
  Vol.~4, 1978, pp. 209--305.

\bibitem{ogden-coupling-09}
R.~Ogden, Incremental elastic motions superimposed on a finite deformation in
  the presence of an electromagnetic field, International Journal of Non-Linear
  Mechanics 44~(5) (2009) 570--580.

\bibitem{saxena-coupledlargedefo-13}
P.~Saxena, On the general governing equations of electromagnetic acoustic
  transducers, Archive of Mechanical Engineering 60~(2) (2013) 231--246.

\bibitem{stiemer-ale-09}
M.~Stiemer, J.~Unger, B.~Svendsen, H.~Blum, {An arbitrary Lagrangian Eulerian
  approach to the three-dimensional simulation of electromagnetic forming},
  Computer Methods in Applied Mechanics and Engineering 198~(17) (2009)
  1535--1547.

\bibitem{abali-largedefo-18-a}
B.~E. Abali, A.~F. Queiruga, Theory and computation of electromagnetic fields
  and thermomechanical structure interaction for systems undergoing large
  deformations, \url{https://arxiv.org/abs/1803.10551}, accessed: 2018-12-20.

\bibitem{ethiraj-coupling-16}
G.~Ethiraj, C.~Miehe, Multiplicative magneto-elasticity of magnetosensitive
  polymers incorporating micromechanically-based network kernels, International
  Journal of Engineering Science 102 (2016) 93--119.

\bibitem{miehe-coupling-16}
C.~Miehe, D.~Vallicotti, S.~Teichtmeister, {Homogenization and multiscale
  stability analysis in finite magneto-electro-elasticity. Application to soft
  matter EE, ME and MEE composites}, Computer Methods in Applied Mechanics and
  Engineering 300 (2016) 294--346.

\bibitem{bayat-coupling-18}
A.~Bayat, F.~Gordaninejad, Characteristic volume element for randomly
  particulate magnetoactive composites, Journal of Engineering Materials and
  Technology 140~(1) (2018) 011003.

\bibitem{homsi-DG-17}
L.~Homsi, L.~Noels, {A discontinuous Galerkin method for non-linear
  electro-thermo-mechanical problems: application to shape memory composite
  materials}, Meccanica (2017) 1--45.

\bibitem{boatti-smp-16}
E.~Boatti, G.~Scalet, F.~Auricchio, {A three-dimensional finite-strain
  phenomenological model for shape-memory polymers: Formulation, numerical
  simulations, and comparison with experimental data}, International Journal of
  Plasticity 83 (2016) 153--177.

\bibitem{wriggers-fem-08}
P.~Wriggers, Nonlinear finite element methods, Springer Science \& Business
  Media, 2008.

\bibitem{penfield-electrodynamics-63}
P.~L. Penfield~Jr, L.~Chu, H.~A. Haus, Electrodynamics of moving media, Tech.
  rep., Research Laboratory of Electronics (RLE) at the Massachusetts Institute
  of Technology (MIT) (1963).

\bibitem{fano-electromagnetism-68}
R.~M. Fano, L.~J. Chu, R.~B. Adler, Electromagnetic fields, energy, and forces,
  MIT Press, 1968.

\bibitem{jackson-electrodynamics-98}
J.~D. Jackson, Classical Electrodynamics, 3rd Edition, Wiley, 1998.

\bibitem{castro-electromagnetism-14}
A.~B. de~Castro, D.~G{\'o}mez, P.~Salgado, Mathematical models and numerical
  simulation in electromagnetism, Vol.~74, Springer, 2014.

\bibitem{simo-fem-06}
J.~C. Simo, T.~J. Hughes, Computational inelasticity, Vol.~7, Springer Science
  \& Business Media, 2006.

\bibitem{belytschko-fem-13}
T.~Belytschko, W.~K. Liu, B.~Moran, K.~Elkhodary, Nonlinear finite elements for
  continua and structures, John Wiley \& Sons, 2013.

\bibitem{hiptmair-mqs-05}
R.~Hiptmair, O.~Sterz, Current and voltage excitations for the eddy current
  model, International Journal of Numerical Modelling: Electronic Networks,
  Devices and Fields 18~(1) (2005) 1--21.

\bibitem{bossavit-cem-98}
A.~Bossavit, Computational Electromagnetism. Variational Formulations,
  Complementarity, Edge Elements, Academic Press, 1998.

\bibitem{scorretti-formulations-12}
R.~Scorretti, R.~V. Sabariego, L.~Morel, C.~Geuzaine, N.~Burais, L.~Nicolas,
  Computation of induced fields into the human body by dual finite element
  formulations, IEEE Transactions on Magnetics 48~(2) (2012) 783--786.

\bibitem{bachinger-cem-05}
F.~Bachinger, U.~Langer, J.~Sch{\"o}berl, Numerical analysis of nonlinear
  multiharmonic eddy current problems, Numerische Mathematik 100~(4) (2005)
  593--616.

\bibitem{hughes-fem-12}
T.~J. Hughes, The finite element method: linear static and dynamic finite
  element analysis, Courier Corporation, 2012.

\bibitem{fish-fem-07}
J.~Fish, T.~Belytschko, A first course in finite elements, John Wiley, 2007.

\bibitem{dular-getdp-98}
P.~Dular, C.~Geuzaine, F.~Henrotte, W.~Legros, A general environment for the
  treatment of discrete problems and its application to the finite element
  method, IEEE Transactions on Magnetics 34~(5) (1998) 3395--3398.

\bibitem{rodriguez-10-eddycurrents}
A.~A. Rodr{\'\i}guez, A.~Valli, Eddy Current Approximation of Maxwell
  Equations: Theory, Algorithms and Applications, Vol.~4, Springer Science \&
  Business Media, 2010.

\bibitem{field-coil}
The electromagnetic fields generated by a long solenoid,
  \url{http://hyperphysics.phy-astr.gsu.edu/hbase/magnetic/solenoid.html},
  accessed: 2018-06-04.

\end{thebibliography}
\end{document}